\documentclass[12pt,leqno]{article}
\usepackage{amsfonts}
\usepackage{amsmath, amsthm, amsfonts, amssymb, color}
\usepackage{mathrsfs}
\usepackage{color}
\newtheorem{thm}{Theorem}[section]
\newtheorem{cor}[thm]{Corollary}

\theoremstyle{definition}

\newcommand{\scr}[1]{\mathscr #1}
\definecolor{wco}{rgb}{0.5,0.2,0.3}

\numberwithin{equation}{section} \theoremstyle{remark}

\newcommand{\ua}{\uparrow}

\title{{\bf    Asymptotic Couplings by Reflection and Applications for  Non-Linear Monotone SPDES  }\footnote{Supported in
 part by  NNSFC(11131003) and Lab. Math. Com. Sys.} }
\author{
{\bf     Feng-Yu Wang  }\\
\footnotesize{ School of Mathematical Sciences,
Beijing Normal
University, Beijing 100875, China}\\
 \footnotesize{ Department of Mathematics,
Swansea University, Singleton Park, SA2 8PP, United Kingdom}\\
\footnotesize{  wangfy@bnu.edu.cn, F.-Y.Wang@swansea.ac.uk}}
\begin{document}
\allowdisplaybreaks
\def\R{\mathbb R}  \def\ff{\frac} \def\ss{\sqrt} \def\B{\mathbf
B}
\def\N{\mathbb N} \def\kk{\kappa} \def\m{{\bf m}}
\def\ee{\varepsilon}\def\ddd{D^*}
\def\dd{\delta} \def\DD{\Delta} \def\vv{\varepsilon} \def\rr{\rho}
\def\<{\langle} \def\>{\rangle} \def\GG{\Gamma} \def\gg{\gamma}
  \def\nn{\nabla} \def\pp{\partial} \def\E{\mathbb E}
\def\d{\text{\rm{d}}} \def\bb{\beta} \def\aa{\alpha} \def\D{\scr D}
  \def\si{\sigma} \def\ess{\text{\rm{ess}}}
\def\beg{\begin} \def\beq{\begin{equation}}  \def\F{\scr F}
\def\Ric{\text{\rm{Ric}}} \def\Hess{\text{\rm{Hess}}}
\def\e{\text{\rm{e}}} \def\ua{\underline a} \def\OO{\Omega}  \def\oo{\omega}
 \def\tt{\tilde} \def\Ric{\text{\rm{Ric}}}
\def\cut{\text{\rm{cut}}} \def\P{\mathbb P} \def\ifn{I_n(f^{\bigotimes n})}
\def\C{\scr C}      \def\aaa{\mathbf{r}}     \def\r{r}
\def\gap{\text{\rm{gap}}} \def\prr{\pi_{{\bf m},\varrho}}  \def\r{\mathbf r}
\def\Z{\mathbb Z} \def\vrr{\varrho} \def\ll{\lambda}
\def\L{\scr L}\def\Tt{\tt} \def\TT{\tt}\def\II{\mathbb I}
\def\i{{\rm in}}\def\Sect{{\rm Sect}}  \def\H{\mathbb H}
\def\M{\scr M}\def\Q{\mathbb Q} \def\texto{\text{o}} \def\LL{\Lambda}
\def\Rank{{\rm Rank}} \def\B{\scr B} \def\i{{\rm i}} \def\HR{\hat{\R}^d}
\def\to{\rightarrow}\def\l{\ell}\def\iint{\int}
\def\EE{\scr E}\def\no{\nonumber}
\def\A{\scr A}\def\V{\mathbb V}\def\osc{{\rm osc}}
\def\BB{\scr B}\def\Ent{{\rm Ent}}

\maketitle

\begin{abstract} Asymptotic couplings by reflection are constructed for a class of non-linear monotone SPDES (stochastic partial differential equations). As applications, the  gradient/H\"older estimates as well as the exponential convergence are  derived for the associated Markov semigroup. The main results are illustrated by stochastic generalized  porous media equations,   stochastic $p$-Laplacian equations, and stochastic generalized fast-diffusion equations. We emphasize that the  gradient estimate is   studied at the first time    for these equations.

\end{abstract} \noindent
 \vskip 2cm

\section{Introduction}

The study of non-linear monotone SPDES (stochastic partial differential equations) goes back to the pioneering work of Pardoux \cite{P1,P2}, see \cite{KR} for an extensive literature on the  existence and uniqueness of solutions as well as It\^o's formula for  the norm of solutions. This type SPDES cover  a number of important models including the stochastic  generalized porous media/fast-diffusion equations and the stochastic $p$-Laplacian equations, which have been intensively investigated in recent years.  Among many other references on this topic we would like to mention \cite{DRRW,RRW,LR1,LR2} for the existence and uniqueness of strong solutions, \cite{BBDR, DR} for  weak solutions, \cite{W07, LW08, L, W14} for the Harnack inequalities and exponential ergodicity, and \cite{BDR, RW13, Gess}   for the extinction properties.
However, so far there is no any  gradient estimates for the associated Markov semigroup.
Below we briefly recall a general formulation of the model and explain the difficulty for the study of gradient estimates.

Let $\V\subset \H\subset \V^*$ be a Gelfand triple, i.e. $(\H,\<\cdot,\cdot\>, \|\cdot\|)$ is a separable Hilbert space, $\V$ is a reflexive Banach space continuously and densely embedded into $\H$, and $\V^*$ is the duality of $\V$ with respect to $\H$.  Let $_{\V^*}\<\cdot,\cdot\>_\V$ be the dualization between   $\V$ and   $\V^*$. We have $_{\V^*}\<u,v\>_\V=\<u,v\>_\H$ for $u\in\H$ and $v\in\V.$ Let $\L_{HS}(\H)$ denote the space of all Hilbert-Schmidt operators on $\H$. Let $W=(W_t)_{t\ge 0}$ be a cylindrical Brownian motion on $\H$,  i.e. $W_t:= \sum_{i=1}^\infty B_t^i e_i$ for an orthonormal basis $\{e_i\}_{i\ge 1}$ of $\H$ and a sequence of independent one-dimensional Brownian motions $\{B_t^i\}_{i\ge 1}$.
 Consider the following stochastic equation:
$$ \d X_t= A(t,X_t)\d t + B(t,X_t) \d W_t,$$
 where
 $A: [0,\infty)\times \V \to \V^*$ and $B: [0,\infty)\times\V\to    \L_{HS}(\H)$ are measurable.    Under some assumptions (see e.g. \cite{LR1}), this equation has a unique solution for any initial point in $\H$. Let $(P_{s,t})_{s\le t}$ be the associated Markov semigroup.  Below we explain the  difficulty for deriving   gradient estimates on $P_{s,t}$ for $t>s$:
 \beg{enumerate}\item[(i)] Since $A$   takes value  in the much larger space $\V^*$   than the state space $\H$, the drift term is highly singular on $\H$, so that the Malliavin calculus does not apply.
 \item[(ii)] Since $B$ is Hilbert-Schmidt, the noise part is too weak to imply gradient estimates as in the semi-linear SPDE case where $QQ^*\ge\ll I$ for some constant $\ll>0$ (see \cite{RW10, WZ}).
   \end{enumerate}

 The purpose of this paper is to investigate   gradient estimates on $P_{s,t}$ using the coupling by reflection. This coupling was introduced by Lindvall and Rogers \cite{LR} and further developed in \cite{CL} and many other papers, and is  optimal  (i.e. the coupling time is minimal) for finite-dimensional Brownian motions.  To investigate gradient estimates for finite-dimensional SDES with multiplicative noise, this coupling was modified in \cite{PW} such that the reflection only occurs to an additive noise part, and the argument has been further developed in \cite{WZ}   for semi-linear SPDES.  However, due to the above two points (i) and (ii), this argument does not apply  to the present non-linear monotone SPDES. Indeed,  since $B$ is Hilbert-Schmidt, the stochastic equation associated to the coupling by reflection is not solvable in the literature (see Section 2 for details). To overcome this difficulty, we will construct a sequence of couplings to approximate the desired $``$coupling by reflection", and we call them asymptotic couplings by reflection.
 Although we can not prove the convergence of these couplings (i.e. the coupling by reflection is not yet well constructed for the present model), the asymptotic couplings by reflection  are already enough to imply gradient/H\"older estimates on the associated Markov semigroup, see Theorem \ref{T1.1} and Theorem \ref{T1.3} below for details.

 In the spirit of  \cite{PW, WZ}, we will split the noise into a multiplicative part and an additive part, and only construct the coupling by reflection for the additive part. To this end, we consider the equation
 \beq\label{E1} \d X_t= A(t,X_t)\d t + B(t, X_t)\d W_t^{(1)} + Q\d W_t^{(2)},\ \ t\in [0,T], \end{equation}
 where $T>0$ is a fixed constant, $Q\in \L_{HS}(\H)$, $$A: [0,T]\times \V\to \V^*, \ \ \  B: [0,T]\times\V\to \L_{HS}(\H)$$ are measurable, and $(W_t^{(1)})_{t\in [0,T]}$ and $(W_t^{(2)})_{t\in [0,T]}$ are two independent cylindrical Brownian motions on $\H$. To construct the asymptotic couplings by reflection, we need the following  assumptions for a fixed constant $r>0$:
   \beg{enumerate} \item[$(A1)$] (Monotonicity) There exists a constant $K\ge 0$   such that for any $t\in [0,T], v_1,v_2\in\V,$
 $$ _{\V^*}\<A(t, v_1)-A(t, v_2), v_1-v_2\>_\V + \ff 1 2  \|B(t, v_1)-B(t,v_2)\|_{HS}^2  \le K_1\|v_1-v_2\|^2.$$
 \item[$(A2)$] (Hemicontinuity) For any  $t\in [0,T]$ and  $v_1,v_2,v\in\V$,\,
 $\R\ni s\mapsto {_{\V^*}\<}A(t, v_1+s v_2), v\>_\V$ is continuous.
  \item[$(A3)$] (Coercivity) There exist constants $C, \theta>0$ such that
 $$_{\V^*}\<A(t, v),v\>_\V +\ff 1 2\|B(t,v)\|_{HS}^2 \le C(1+\|v\|^2)-\theta\|v\|_\V^{r+1},\ \    t\in [0,T], v\in \V.$$
 \item[$(A4)$] (Growth) There exists a constant $c>0$ such that
 $$ |_{\V^*}\<A(t,v),u\>_\V|\le c \big(1+\|v\|_\V^{r+1}+\|u\|_\V^{r+1}+\|u\|^2+\|v\|^2\big),\ \ t\in [0,T], u,v\in\V.$$ \end{enumerate}
Then \eqref{E1} has a unique solution with $X_0=x$   for any initial point $x\in\H$,  and we denote the solution by $(X_t(x))_{t\in [0,T]}$ (see e.g. \cite{LR1}). Recall that a continuous $\H$-valued adapted process $X:=(X_t)_{t\in[0,T]}$ is called a   solution to (\ref{E1}), if
$$ \E\int_0^T \big( \|X_t\|^2+ \|X_t\|_\V^{r+1}\big)\,\d t  <\infty,\ \ T>0,$$   and $\P$-a.s.
$$X_t=X(0) +\int_0^t A(s,X_s)\,\d s +\int_0^t B(s,X_s)\d W_s^{(1)}+  Q   W_t^{(2)},\ \ t\in [0,T],$$ where the Bochner integral $\int_0^t A(s,X_s)\,\d s$  is defined on $\V^*$ but it  takes values in $\H$ for all $t\in [0,T]$, and  in this integral as well as in the stochastic integral $\int_0^t B(x,X_s)\d W_s^{(1)}$ we have identified $X$ and its proper modification $\bar X\in L^2([0,T]\times \OO\to \H)\cap L^{1+r}([0,T]\times\OO\to\V)$ used in \cite{LR1}.

\

Let $(P_{s,t})_{T\ge t\ge s\ge 0}$ be the associated Markov semigroup, i.e.
$$P_{s,t} f(x):= \E f(X_{s,t}(x)),\ \ f\in \B_b(\H), T\ge t\ge s\ge 0, x\in\H,$$ where $(X_{s,t}(x))_{t\in [s,T]}$ solves the equation \eqref{E1}
from time $s$ with $X_{s,s}(x)=x.$ By the uniqueness of solutions, we have the semigroup property
$$P_{s,t}= P_{s, t'}P_{t',t},\ \ T\ge t\ge t'\ge s\ge 0.$$  We simply denote $P_t=P_{0,t}, t\in [0,T].$

To derive gradient estimates on $P_{s,t}$, we also need the intrinsic norm induced by $Q$:
$$\|u\|_Q= \|(QQ^*)^{-\ff 1 2}u\|:=\inf\big\{\|x\|: x\in \H, (QQ^*)^{\ff 1 2} x=u\big\},\ \ u\in\H,$$ where we set $\inf\emptyset =\infty$ by convention.

\

The remainder of the paper is organized as follows. In Section 2, we introduce the main results of the paper and compare them with some known ones. In Section 3, we construct the asymptotic couplings by reflection and use them to estimate  $|P_tf(x)-P_tf(y)|$ for $f\in \B_b(\H), t>0$ and $x,y\in\H$.
In Sections 4 and 5, we prove our main results using asymptotic couplings by reflections for  $r\ge 1$ and $r\in (0,1)$ respectively. Finally, in Section 6, some specific examples are presented to illustrate the main results.

\section{Main results}

According to \cite{L, LW08, W07,Wbook} where the Harnack inequalities for monotone SPDES are studied, we shall need the following stronger versions of $(A1)$ for $r\ge 1$ and $r\in (0,1)$ respectively, which are easy to check in applications (see Section 6 for details).

\beg{enumerate}\item[$(A1')$]    $r\ge 1$, and  there exist  constants $K,\theta>0$ and $\kk>r-1$ such that
 \beg{equation*}\beg{split} &  _{\V^*}\<A(t, v_1)-A(t, v_2), v_1-v_2\>_\V + \ff 1 2  \|B(t, v_1)-B(t,v_2)\|_{HS}^2\\
   &\le K|v_1-v_2\|^2-\theta\|v_1-v_2\|^{r+1-\kk}\|v_1-v_2\|_Q^\kk,\ \  t\in [0,T], v_1,v_2\in\V.\end{split}\end{equation*}
 \item[$(A1'')$]    $r\in (0, 1)$, and  there exist  constants $K,\theta, \kk>0$ such that
 \beg{equation*}\beg{split} & _{\V^*}\<A(t, v_1)-A(t, v_2), v_1-v_2\>_\V + \ff 1 2  \|B(t, v_1)-B(t,v_2)\|_{HS}^2\\
   &\le K\|v_1-v_2\|^2-\theta\ff{\|v_1-v_2\|^{2-\kk}\|v_1-v_2\|_Q^\kk}{(\|v_1\|_\V\lor \|v_2\|_\V)^{1-r}},\ \  t\in [0,T], v_1,v_2\in\V.\end{split}\end{equation*}\end{enumerate}

\subsection{For $r\ge 1$}

\beg{thm}\label{T1.1} Assume $(A1')$ and $(A2)$-$(A4)$.
\begin{enumerate} \item[$(1)$] If $\kk>2r $, then there exists a constant $C>0$ such that for any $0\le s<t\le T$ and $f\in \B_b(\H)$,
$$\|\nn P_{s,t}f\|(x):=\limsup_{y\to x} \ff{|P_{s,t} f(x)-P_{s,t}f(y)|}{\|x-y\|}
\le \ff{C\|f\|_\infty}{\{(t-s)\land 1\}^{\ff{\kk+2}{2\kk+2-2r}}},\ \ x\in\H.$$
\item[$(2)$] If $\kk=2r,$ then there exists a constant $C>0$ such that   for any $0\le s<t\le T$ and $f\in \B_b(\H)$,
$$|P_{s,t}f(x)-P_{s,t}f(y)|\le \ff{C\|f\|_\infty}{(t-s)\land 1} \|x-y\|\big\{\log(\|x-y\|^{-1}+\e)\big\}^{\ff r{r+1}},\ \   x,y\in \H.$$
\item[$(3)$] If $\kk\in (r-1, 2r)$, then there exists a constant $C>0$ such that for any $0\le s<t\le T$ and $f\in \B_b(\H)$,
$$|P_{s,t}f(x)-P_{s,t}f(y)|\le \ff{C\|f\|_\infty\|x-y\|^{\ff{2(\kk+1-r)}{\kk+2}}}{(t-s)\land 1},\ \  x,y\in \H.$$
\end{enumerate}
\end{thm}

 \paragraph{Remark 2.1.} (1) When $\H=\R^d$ is finite-dimensional, $P_{t}$ becomes an elliptic diffusion semigroup. In this case it is well known that the sharp short time behavior of $|\nn P_t|$ is of order $\ff 1 {\ss t}.$  Since in this case $(A1')$   holds for any $\kk>r-1$, the gradient estimate in Theorem \ref{T1.1}(1) is sharp for short time when $\kk\to\infty$.

 (2) Using coupling by change of measures constructed in \cite{W07}, it is shown in \cite[Theorem 2.2.1]{Wbook} that $(A1')$ and $(A2)$-$(A4)$ with $\kk\in (r-1,\infty)\cap [2,\infty)$ imply the log-Harnack inequality
 $$P_t \log f(y)\le \log P_t f(y) +\ff{C\|x-y\|^{\ff{2(\kk+1-r)}\kk}}{(t\land 1)^{\ff{\kk+2}{\kk}}}$$
 for all positive $f\in \B_b(\H), t\in (0,T]$ and $x,y\in\H.$ According to \cite[Proposition 2.3]{ATW}, this implies
 $$|P_tf(x)-P_tf(y)|\le \ff{C\|f\|_\infty\|x-y\|^{\ff{\kk+1-r}\kk}}{(t\land 1)^{\ff{\kk+2}{2\kk}}},\ \  x,y\in \H.$$
Since when $\kk> 2$ we have $\ff{\kk+1-r}\kk< \ff{2(\kk+1-r)}{\kk+2}$,  this  H\"older continuity is worse than that in Theorem \ref{T1.1}(3).  Moreover, when
$\kk\ge 2r$, results in  Theorem \ref{T1.1}(1)-(2) are much stronger than the H\"older continuity.

\

As an application of Theorem \ref{T1.1}, we have the following result on the exponential convergence of $P_t$. When $B=0$ and $A(t,v)$ does not depend on $t$, this property has been investigated in \cite{L, W14}.

\beg{thm}\label{T1.2} Assume   $(A1')$ and $(A2)$-$(A4)$  for all $T>0$. \beg{enumerate}
\item[$(1)$] If $r>1$,    then there exists constants $C,\ll>0$ such that
\beq\label{UE} \sup_{x,y\in \H} |P_t f(x)-P_t f(y)| \le C\|f\|_\infty \e^{-\ll t},\ \ t\ge 0, f\in \B_b(\H).\end{equation}
\item[$(2)$] If $r=1$ and there exists a constant $K>0$ such that
\beq\label{EC}\beg{split} & _{\V^*}\<A(t,v_1)-A(t,v_2), v_1-v_2\>_\V +\ff 1 2 \|B(t,v_1)-B(t,v_2)\|_{HS}^2\\
 &\le -K\|v_1-v_2\|^2,\ \ t\ge 0, v_1,v_2\in\V,\end{split}\end{equation} then there exists a constant $C>0$   such that
\beq\label{EX}  \ff{|P_tf(x)-P_t f(y)|}{C \|f\|_\infty}\le \beg{cases}   \|x-y\|\e^{-Kt}, &\text{if}\ \kk>2,\\
  \|x-y\|\e^{-Kt}\ss{\log(\e + \ff{\e^{Kt}}{\|x-y\|})}, &\text{if}\ \kk=2,\\
 \|x-y\|^{\ff {2\kk}{\kk+2}} \e^{-\ff{2K\kk t}{\kk+2}}, &\text{if}\ \kk\in (0,2)\end{cases} \end{equation} holds for all $x,y\in\H, t\ge 0$ and $0\ne f\in \B_b(\H).$ \end{enumerate} \end{thm}

\paragraph{Remark 2.2.} Let $P_t(x,\cdot)$ be the distribution of $X_t(x)$. Then \eqref{UE} is equivalent to
$$\sup_{x,y\in \H} \|P_t(x,\cdot)-P_t(y,\cdot)\|_{var} \le C\e^{-\ll t},\ \ t\ge 0,$$ where $\|\cdot\|_{var}$ is the total variational norm.
If $A(t,v)$ and $B(t,v)$ are independent of $t$, then in the situation of Theorem  \ref{T1.2}(1),  $P_t$ has a unique invariant probability $\mu$ such that
$$ \sup_{x\in\H} \|P_t(x,\cdot) -\mu\|_{var}\le C\e^{-\ll t},\ \ t\ge 0,$$ which is known as the strong ergodicity.
Moreover, $(A3)$  implies   $\mu(\|\cdot\|^2)<\infty$, so that when $r=1$ Theorem \ref{T1.2}(2) implies
$$\mu\big((P_t f-\mu(f))^2\big)\le C \|f\|_\infty^2 (1+ t 1_{\{\kk=2\}})\exp\bigg[-\ff{2Kt}{1\lor \ff{2+\kk}{2\kk}}\bigg],\ \ t\ge 0, f\in \B_b(\H).$$ Therefore, Theorem \ref{T1.2} generalizes,  and improves when $r=1$,  assertions in \cite[Theorem 1.5]{L} where $B=0$ and $A(t,v)=A(v)$ is considered.

\subsection{For $r\in (0,1)$}
We call $B$ bounded, if
$$\|B\|_\infty:= \sup\{\|B(t,x)z\|:\ t\in [0,T], x\in\V, z\in\H, \|z\|\le 1\}<\infty.$$

\beg{thm}\label{T1.3} Assume   $(A1'')$ and $(A2)$-$(A4)$.
\begin{enumerate} \item[$(1)$] If $\kk>2$, then for any $p>0$ there exists a constant $C_p>0$ such that
$$ | P_{s,t}f (x)- P_{s,t}f (y)|
\le \ff{C_p\|f\|_\infty\|x-y\|^{\ff{2p\kk}{2p\kk+\kk+2}}(1+\|x\|^2+\|y\|^2)^{\ff{2p(1-r)}{(1+r)(2p\kk+\kk+2)}}}{\{(t-s)\land 1\}^{\ff{p(4+\kk +\kk r)}{(r+1)(2p\kk+\kk+2)}}}$$ holds for  for any $0\le s<t\le T,$  $f\in \B_b(\H)$ and $x,y\in\H.$ If moreover $B$ is bounded, then there exists a constant $C>0$ such that
$$ | P_{s,t}f (x)- P_{s,t}f (y)|
\le \ff{C\|f\|_\infty\|x-y\|\{1+\|x\|^2+\|y\|^2+\log (\e+\|x-y\|^{-1})\}^{\ff{1-r}{\kk(1+r)}}}{\{(t-s)\land 1\}^{\ff{4+\kk +\kk r}{2\kk (r+1)}}}$$ holds for  for any $0\le s<t\le T,$  $f\in \B_b(\H)$ and $x,y\in\H.$
\item[$(2)$] If $\kk=2,$   then for any $p>0$ there exists a constant $C_p>0$ such that
\beg{equation*}\beg{split} &| P_{s,t}f (x)- P_{s,t}f (y)|\\
&\le \ff{C_p\|f\|_\infty\big\{\|x-y\|^2 \log(\e+\|x-y\|^{-1})\big\}^{\ff{p}{2(p+1)}}  (1+\|x\|^2+\|y\|^2)^{\ff{p(1-r)}{2(p+1)(1+r)}}}{\{(t-s)\land 1\}^{\ff{p(3+ r)}{2(p+1)(1+r)}}} \end{split}\end{equation*} holds for  for any $0\le s<t\le T,$  $f\in \B_b(\H)$ and $x,y\in\H.$ If moreover $B$ is bounded, then  there exists a constant $C>0$ such that
\beg{equation*}\beg{split} &| P_{s,t}f (x)- P_{s,t}f (y)|\\
&\le \ff{C\|f\|_\infty\|x-y\|\ss{\log (\e+\|x-y\|^{-1})}\big\{1+\|x\|^2+\|y\|^2+\log (\e+\|x-y\|^{-1})\big\}^{\ff{1-r}{2(1+r)}}}{\{(t-s)\land 1\}^{\ff{3+r}{2 (r+1)}}}\end{split}\end{equation*} holds for  for any $0\le s<t\le T,$  $f\in \B_b(\H)$ and $x,y\in\H.$
\item[$(3)$] If $\kk\in (0, 2)$, then for any $p>0$ there exists a constant $C_p>0$ such that  for any $0\le s<t\le T$ and $f\in \B_b(\H)$,
$$ | P_{s,t}f (x)- P_{s,t}f (y)|
\le \ff{C_p\|f\|_\infty\|x-y\|^{\ff{2p\kk}{(p+1)(\kk+2)}}(1+\|x\|^2+\|y\|^2)^{\ff{2p(1-r)}{(1+r)(2p\kk+\kk+2)}}}{\{(t-s)\land 1\}^{\ff{p(4+\kk +\kk r)}{(r+1)(2p\kk+\kk+2)}}},\ \ x,y\in\H.$$ If moreover $B$ is bounded, then there exists a constant $C>0$ such that  for any $0\le s<t\le T$ and $f\in \B_b(\H)$,
$$ | P_{s,t}f (x)- P_{s,t}f (y)|
\le \ff{C\|f\|_\infty\|x-y\|^{\ff{2\kk}{\kk+2}}\{1+\|x\|^2+\|y\|^2+\log (\e+\|x-y\|^{-1})\}^{\ff{1-r}{\kk(1+r)}}}{\{(t-s)\land 1\}^{\ff{4+\kk +\kk r}{2\kk (r+1)}}}$$ holds for all $x,y\in\H.$
\end{enumerate}
\end{thm}

If the constant $C_p$ in Theorem \ref{T1.3}(1) is bounded as  $p\to\infty$, then ny letting $p\to\infty$ we obtain the gradient estimate
$$\|\nn P_{s,t}f\|(x)
\le \ff{C\|f\|_\infty(1+\|x\|^2)^{\ff{(1-r)}{\kk(1+r)}}}{\{(t-s)\land 1\}^{\ff{4+\kk(r+1)}{2\kk(r+1)}}},\ \ f\in\B_b(\H), x\in\H$$ for some constant $C>0.$
 However, in this paper we can not prove the boundedness of $(C_p)_{p>0}$.  Nevertheless, the following result implies this gradient estimate when $B=0$ and $\kk\ge \ff 4 {r+1}$. The proof of this result is due to the log-Harnack inequality presented in \cite[Theorem 2.3.1]{Wbook}. But the gradient estimate for non-constant $B$ or $\kk <\ff 4{r+1}$ is still unknown.

 \beg{prp} \label{P*} Assume $(A1'')$ and $(A2)$-$(A4)$. If    $B=0$ and   $\kk\ge \ff 4 {r+1}$, then
 \beq\label{NNB} \|\nn P_{s,t}f\|^2(x)
\le \ff{C(P_{s,t} f^2(x)) (1+\|x\|^2)^{\ff{2(1-r)}{\kk(1+r)}}}{\{(t-s)\land 1\}^{\ff{4+\kk(r+1)}{\kk(r+1)}}},\ \ f\in \B_b(\H), x\in\H\end{equation}
 holds for some constant $C>0$. \end{prp}

 \beg{proof} According to \cite[Theorem 2.3.1]{Wbook},   $(A1'')$ and $(A2)$-$(A4)$ for $B=0$ and   $\kk\ge \ff 4 {r+1}$ imply the log-Harnack inequality
$$ P_t\log f(y)\le \log P_t f(x) +\ff{C \|x-y\|^2 (1+\|x\|^2+\|y\|^2)^{\ff{2(1-r)}{\kk(1+r)}}}{(t\land 1)^{\ff{\kk(1+r)+4}{\kk(1+r)}}},\ \ t\in (0,T], x,y\in \H$$   for some constant $C>0$ and all positive $f\in \B_b(\H).$  Using the argument in the proof of \cite[Proposition 2.3]{ATW}, it is easy to see that this implies the desired gradient estimate \eqref{NNB}.  \end{proof}

\section{Asymptotic couplings by reflection} Throughout this section, we assume $(A1)$-$(A4)$ and let $Q$ be symmetric with $Q\ge 0$ and Ker$(Q)=\{0\}.$  According to \cite{PW, WZ}, we shall construct the coupling by reflection for the additive noise part. To this end, let
$$\si(u,v)= \ff{\{Q^{-1}(u-v)\}\otimes \{Q^{-1}(u-v)\}}{\|u-v\|_Q^2},\ \ \ u,v\in  Q(\H).$$
Then the stochastic equation for the coupling process $(X_t,Y_t)$ is formulated as follows:
$$\beg{cases} \d X_t = A(t,X_t)\d t + B(t,X_t) \d W_t^{(1)} +Q \d W_t^{(2)}, \\
\ \\
\d Y_t = A(t,Y_t)\d t + B(t,Y_t) \d W_t^{(1)} +Q\big(I-2 \si(X_t,Y_t)\big) \d W_t^{(2)},\end{cases}$$ where in the second equation   $W_t^{(2)}$ is reflected along the direction  $Q^{-1}(X_t-Y_t).$

However, since $Q$ is Hilbert-Schmidt, $Q\si(\cdot): \H^2\to \L_{HS}(\H)$ is  not a well-defined continuous map, so that existing results on the existence and uniqueness of solutions for monotone SPDES do not apply, see \cite{LR1, LR2, RRW}. To overcome this difficulty, we use $(Q+\ff 1 n I)^{-1}$ to replace $Q^{-1}$: for any $n\ge 1$, let
$$\si_n(u,v)= \ff{\{(Q+\ff 1 n I)^{-1}(u-v)\}\otimes \{(Q+\ff 1 n I)^{-1}(u-v)\}}{\|(Q+\ff 1 n I)^{-1}(u-v)\|^2},\ \ t\in [0,T], u,v\in\H.$$
Since $Q$ is non-negative definite, for any $n\ge 1$ the norm $\|(Q+\ff 1 n I)^{-1} x\|$ is equivalent to $\|x\|$. Moreover, to manage the singularity of $\si_n(u,v)$ on the diagonal $\{(u,u): u\in\H\}$, we introduce a cut-off function
\beq\label{H} h(s)=\beg{cases} 0, & s\in [0,\ff 1 2],\\
1-\exp[-(r-\ff 1 2)/(r-1)], & s\in (\ff 1 2,1),\\
1, &s\ge 1.\end{cases}\end{equation}
Obviously, $h, \ss{1-h^2}\in C_b^1([0,\infty))$.

Now, for any $n\ge 1$, consider the coupling stochastic equation
\beq\label{EN} \beg{cases} \d X_t^n = A(t,X_t^n)\d t + B(t,X_t^n) \d W_t^{(1)} +Q\ss{1-h(n\|X_t^n-Y_t^n\|)^2}\,\d W_t^{(2)}
 \\ \qquad\qquad
 +Q  h(n\|X_t^n-Y_t^n\|)\d W_t^{(3)}, \\
 \ \\
\d Y_t^n = A(t,Y_t^n)\d t + B(t,Y_t^n) \d W_t^{(1)} +Q\ss{1-h(n\|X_t^n-Y_t^n\|)^2}\,\d W_t^{(2)}
  \\
  \qquad\qquad +Q h(n\|X_t^n-Y_t^n\|)\big(I-2 \si_n(X_t,Y_t)\big) \d W_t^{(3)},\end{cases}\end{equation} where $W_t^{(1)}, W_t^{(2)}$ and $W_t^{(3)}$ are independent cylindrical Brownian motions on $\H$. Since $h(n s)=0$ for $s\le \ff 1{2n}$,   the reflection   occurs only when
  $\|X_t^n-Y_t^n\|>\ff 1 {2n}.$   To see that this equation has a unique solution for any initial point in $\H^2$, for any
  $(u,v)\in \V^2, (v_1,v_2,v_3)\in\H^3$ and $t\in [0,T]$, we let
\beg{equation*}\beg{split} &\bar A(t,(u,v))= (A(t,u), A(t, v)),\\
&\bar B(t, (u,v))(v_1,v_2,v_3)= \Big(B(t,u)v_1+ h(n\|u-v\|)Qv_2+ \ss{1-h(n\|u-v\|)^2}\, Qv_3,\\
&\qquad\qquad\qquad B(t,v)v_1+ h(n\|u-v\|)Q\big(I-2\si_n(u,v)\big)v_2+ \ss{1-h(n\|u-v\|)^2}\, Qv_3\Big).\end{split}\end{equation*} Moreover, let $\bar W_t= (W_t^{(1)}, W_t^{(2)}, W_t^{(3)})$ which is a cylindrical Brownian motion on $\H^3$. Then \eqref{EN} can be reformulated as the following equation on $\bar\H:=\H^2$:
$$\d \bar X_t = \bar A(t, \bar X_t) \d t + \bar B(t, \bar X_t)\d \bar W_t.$$
Let $\bar \V=\V^2, \bar\V^*=(\V^*)^2$. Then $\bar\V\subset\bar\H\subset\bar\V^*$ is a Gelfand triple. Moreover, it is easy to  see from $(A1)$-$(A4)$ and the construction of $\si_n$ that $\bar A$ and $\bar B$ satisfy the following conditions for some constants $K,\theta>0$ depending on $n$, and all $t\in [0,T], \bar v, \bar v_1,\bar v_2\in \bar V$:
\beg{enumerate} \item[{\bf (a)}] (Hemicontinuity)
 $\R\ni s\mapsto {_{\bar\V^*}\<}\bar A(t, \bar v_1+s \bar v_2), v\>_{\bar\V}$ is continuous.
 \item[{\bf (b)}] (Monotonicity)
 $_{\bar\V^*}\<\bar A(t, \bar v_1)-\bar A(t, \bar v_2), \bar v_1-\bar v_2\>_{\bar\V} + \ff 1 2 \|\bar B(t,\bar v_1)-\bar B(t,\bar v_2)\|_{\L_{HS}(\H^3\to \bar\H)}^2 \le K \|\bar v_1-\bar v_2\|_{\bar\H}^2.$
 \item[{\bf (c)}] (Coercivity)
 $_{\bar\V^*}\<\bar A(t, \bar v),\bar v\>_{\bar\V} +\ff 1 2\|\bar B(t,\bar v)\|_{\L_{HS}(\H^3\to\bar\H)}^2 \le K (1+\|\bar v\|_{\bar\H}^2)-\theta \|v\|_{\bar\V}^{r+1}.$
 \item[{\bf (d)}] (Growth)
 $ |_{\bar\V^*}\<\bar A(t,\bar v_1),\bar v_2\>_{\bar\V}|\le K  \big(1+\|\bar v_1\|_{\bar\V}^{r+1}+\|\bar v_2\|_{\bar\V}^{r+1}+\|\bar v_1\|_{\bar\H}^2+\|\bar v_2\|_{\bar\H}^2\big).$ \end{enumerate} Therefore, for any initial point $(x,y)\in \H^2$, the equation \eqref{EN} has a unique solution $(X_t^n(x,y), Y_t^n(x,y))$ staring at $(x,y),$
 see e.g. \cite[Theorem 2.1]{RRW}.
Below we prove that   the solution is a coupling of $X_t(x)$ and $X_t(y)$, i.e. the law of $(X_t^n(x,y))_{t\in [0,T]}$ coincides with that of $(X_t(x))_{t\in [0,T]}$, and the same is true for $(Y^n_t(x,y))_{t\in [0,T]}$ and $(X_t(y))_{t\in [0,T]}.$ Moreover, when $n\to\infty$, these couplings provide an upper bound estimate of $|P_t f(x)-P_tf(y)|.$

\beg{prp}\label{P2.1} Assume $(A1)$-$(A4)$. Let $Q\ge 0$ be symmetric.
\beg{enumerate} \item[$(1)$]   $(X_t^n(x,y), Y_t^n(x,y))_{t\in [0,T]}$ is a coupling of $(X_t(x))_{t\in [0,T]}$ and $(X_t(y))_{t\in [0,T]}$.
\item[$(2)$] Let $T_n^{x,y}= \inf\{t\in [0,T]: X_t^n(x,y)=Y_t^n(x,y)\}$ be  the coupling time, where $\inf \emptyset =\infty$ by convention. Then
 $X_t^n(x,y)= Y_t^n(x,y)$ holds for $t\in [T_n^{x,y},T]$.
\item[$(3)$] Let $\tau_n^{x,y}= \inf\{t\in [0,T]: \|X_t^n(x,y)-Y_t^n(x,y)\|\le \ff 1 n\}.$ Then
$$|P_t f(x)-P_t f(y)|\le \osc(f) \liminf_{n\to\infty} \P(\tau_n^{x,y}>t),\ \ t\in (0,T], f\in \B_b(\H), x,y\in \H,$$ where $\osc(f):=\sup f-\inf f.$
\end{enumerate}
\end{prp}

\beg{proof}  We simply denote $(X_t^n,Y_t^n)= (X_t^n(x,y),Y_t^n(x,y))$,  $X=(X_t)_{t\in [0,T]}$ for a process
$X_t$ on $\H$, $T_n=T_n^{x,y}$ and $\tau_n=\tau_n^{x,y}$.

(1) It is easy to see that
$$\tt W_t^{(2)} := \int_0^t\Big\{\ss{1-h(n\|X_s^n-Y_s^n\|)^2}\,\d W_s^{(2)}
  +  h(n\|X_s^n-Y_s^n\|)\d W_s^{(3)}\Big\}$$ is a cylindrical Brownian motion on $\H$, which is independent of $W^{(1)}$ since, for any $u,v\in\H$,
  the processes $\<u,W_t^{(1)}\>$ and $\<v, \tt W_t^{(2)}\>$  have zero covariation. By \eqref{EN} we have
  $$\d X_t^n= A(t, X_t^n)\d t +B(t, X_t^n)\d W_t^{(1)}+ Q\d\tt W_t^{(2)},\ \ X_0^n=X_0=x.$$ So,   the uniqueness of the weak solutions to \eqref{E1} implies that $X(x)$ and $X^n$ have a common   distribution.

Similarly, since $\si_n$ is symmetric with $(I-2\si_n)^2=I$, the above argument also applies to $Y^n$ and $X(y)$ so that they have a common distribution as well.

(2) Let
$T_n'= \inf\big\{t\in [T_n,T]: \|X_t^n-Y_t^n\|\ge \ff 1{2n}\Big\}.$
Since $h(ns)=0$ for $s\in [0,\ff 1 {2n}]$,   \eqref{EN} yields
$$\d(X_t^n-Y_t^n)= (A(t,X_t^n)-A(t,Y_t^n))\d t+ (B(t,X_t^n)-B(t,Y_t^n))\d W_t^{(1)},\  t\in [T_n,T_n'\land T].$$
By $(A1)$ and It\^o's formula (see e.g. \cite[Theorem A.2]{RRW}), there exists a constant $c>0$ such that
$$\d \|X_t^n-Y_t^n\|^2 \le c\|X_t^n-Y_t^n\|^2 \d t + 2\<(B(t,X_t^n)-B(t,Y_t^n))\d W_t^{(1)}, X_t^n-Y_t^n\>$$ holds for $t\in [T_n,T_n'\land T].$
Thus,
\beg{equation*}\beg{split} &1_{\{T_n\le T\}}\E \big(\|X_{(T_n+t)\land T_n'\land T}^n-Y_{(T_n+t)\land T_n'\land T}^n\|^2\e^{-c ((T_n+t)\land T_n'\land T)}\big|\F_{T_n\land T}\big)\\
&\le 1_{\{  T_n\le T\}}   \E \big(\|X_{T_n}-Y_{T_n}\|^2\e^{-cT_n}\big|\F_{T_n\land T}\big)=0,\ \ t\in [0,T].\end{split}\end{equation*}
This implies  $T_n'=\infty$ and $X_t^n=Y_t^n$   for $t\in [T_n,T].$

(3) Without loss of generality, we may and do assume that $f$ is Lipschitz continuous.   Let $Q_n=Q+\ff 1 n I$.
By \eqref{EN} we have
\beg{equation*}\beg{split} \d(X_t^n-Y_t^n)=&\ (A(t,X_t^n)-A(t,Y_t^n))\d t+ (B(t,X_t^n)-B(t,Y_t^n))\d W_t^{(1)} \\
&\ + 2  h(n\|X_t^n-Y_t^n\|) Q\si_n(X_t^n,Y_t^n) \d W_t^{(3)},\ \ t\in [0,T].
\end{split}\end{equation*} By It\^o's formula (see e.g. \cite[Theorem A.2]{RRW}), this implies
\beg{equation}\label{W*F}\beg{split} &\d \|X_t^n-Y_t^n\|^2 =  \Big\{2_{\V^*}\<A(t,X_t^n)-A(t,Y_t^n), X_t^n-Y_t^n\>+\|B(t,X_t)-B(t,Y_t)\|_{HS}^2 \Big\}\d t\\
& \qquad + \ff{4h(n\|X_t^n-Y_t^n\|)^2\|QQ_n^{-1} (X_t^n-Y_t^n)\|^2}{\|Q_n^{-1} (X_t^n-Y_t^n)\|^2}\,\d t\\
& \qquad + 2\Big\<(B(t,X_t^n)-B(t,Y_t^n))\d W_t^{(1)}, X_t^n-Y_t^n\Big\> \\
&\qquad + \ff{4 h(n\|X_t^n-Y_t^n\|) \<QQ_n^{-1}(X_t^n-Y_t^n), X_t^n-Y_t^n\>}
{\|Q_n^{-1} (X_t^n-Y_t^n)\|^2} \Big\<Q_n^{-1} (X_t^n-Y_t^n), \d W_t^{(3)}\Big\> \end{split}\end{equation}for $t\in [0,T].$
So, by $(A1)$, there exists a constant  $K >0$ such that
\beq\label{D1} \d\|X_t^n-Y_t^n\|\le \Big\{K\|X_t^n-Y_t^n\|  +I_n(X_t^n-Y_t^n)\Big\}\d t
+\d M_t^n,\ t <T_n\land T,\end{equation}where
\beg{equation}\beg{split}\label{D2} &\d M_t^n = \Big\<(B(t,X_t^n)-B(t,Y_t^n))\d W_t^{(1)},\ff{X_t^n-Y_t^n}{\|X_t^n-Y_t^n\|}\Big\> \\
& + \ff{2 h(n\|X_t^n-Y_t^n\|) \<QQ_n^{-1}(X_t^n-Y_t^n), X_t^n-Y_t^n\>}
{\|Q_n^{-1} (X_t^n-Y_t^n)\|\cdot \|X_t^n-Y_t^n\|} \Big\<\ff{Q_n^{-1} (X_t^n-Y_t^n)}{\|Q_n^{-1}(X_t^n-Y_t^n)\|}, \d W_t^{(3)}\Big\>,\end{split}\end{equation} and
$$I_n(v)= \ff{2 h(n\|v\|)^2 }{\|v\|\cdot\|Q_n^{-1}v\|^2}\Big(\|QQ_n^{-1}v\|^2- \ff{\<QQ_n^{-1}v, v\>^2}{\|v\|^2}\Big),\ \ 0\ne v\in\H.$$
Since $Q_n=Q+\ff 1 n I$, we have
\beg{equation*}\beg{split} &\|QQ_n^{-1}v\|^2- \ff{\<QQ_n^{-1}v, v\>^2}{\|v\|^2}= \Big\|v-\ff 1 n Q_n^{-1} v\Big\|^2 - \ff {(\|v\|^2 -\ff 1 n \<Q_n^{-1}v,v\>)^2} {\|v\|^2}\\
&=\ff 1 {n^2}\Big(\|Q_n^{-1}v\|^2 -\ff{\<Q_n^{-1} v,v\>^2}{\|v\|^2}\Big) \le \ff 1{n^2} \|Q_n^{-1}v\|^2.\end{split}\end{equation*}  So,
$$I_n(X_t^n-Y_t^n) \le \ff{2h(n\|X_t^n-Y_t^n\|)^2}{n^2\|X_t^n-Y_t^n\|} \le 2 \|h'\|_\infty^2\|X_t^n-Y_t^n\|,\ \ t< T_n\land T.$$Combining this with \eqref{D1}, we obtain
\beq\label{D3} \d\|X_t^n-Y_t^n\|\le  K'\|X_t^n-Y_t^n\|  \d t
+\d M_t^n,\ \ t <T_n\land T\end{equation} for some constants $K'>0.$ Since by (2) we have $X_t^n=Y_t^n$ for $t\in [T_n,T]$, this implies that
$(\|X_t^n-Y_t^n\|\e^{-K't})_{t\in [0,T]}$ is a supermartingale, so that
\beq\label{*W}\E\big(\e^{-K't}\|X_t^n-Y_t^n\|1_{\{t\ge\tau_n\}}\big)\le \E\big(\|X_{\tau_n}^n- Y_{\tau_n}^n\| 1_{\{t\ge \tau_n\}}\big) \le\ff 1 n,\ \ t\in [0,T].\end{equation}
Combining this with (1) we conclude that  for any Lipschitz function $f$ on $\H$ with Lipschitz constant $L(f)$,
\beg{equation*}\beg{split} &|P_tf(x)-P_t f(y)|=|\E(f(X_t^n)-f(Y_t^n))|\\
&\le \E\big|(f(X_t^n)-f(Y_t^n))1_{\{\tau_n>t\}}\big| + L(f) \e^{K't} \E\big(\e^{-K't}\|X_t^n-Y_t^n\|1_{\{t\ge\tau_n\}}\big)\\
&\le \osc(f) \P(\tau_n>t) + \ff{L(f)\e^{K' t}}n,\ \ \ n\ge 1.\end{split}\end{equation*}
Letting $n\to\infty$ we prove (3).
\end{proof}

\section{Proof of Theorem \ref{T1.1} and Theorem \ref{T1.2}}
Since $\d \tt W^{(2)}_t := (QQ^*)^{-\ff 1 2}Q\d W_t^{(2)}$ gives rise  to a cylindrical Brownian motion on $\H$ independent of $\d W^{(1)}_t,$ and   since $Q\d W_t^{(2)} = \ss{QQ^*}\, \d \tt W_t^{(2)},$ in \eqref{E1}  we may and do assume that $Q$ is symmetric with $Q\ge 0.$

Let
$x,y\in \H$ with $x\ne y$. We simply denote $(X_t^n,Y_t^n)$ the solution to \eqref{EN} for $(X_0^n,Y_0^n)=(x,y)$,  and let $\tau_n=\tau_n^{x,y}.$
Then, according to Proposition \ref{P2.1}(3), the key point for the proof of Theorem \ref{T1.1} is to estimate $\P(\tau_n>t),$  for which we  follow the line of  \cite{CL}. For any $\dd>0,$   let
$$\tau_{n,\dd}=\inf\{t\in [0,T]: \|X_t^n-Y_t^n\|\ge\dd\}.$$ We have  $\tau_{n,\dd}=0$ if $\dd\le \|x-y\|.$ According to \cite{CL},    we  need to estimate $\E(\tau_n\land\tau_{n,\dd}\land t)$ and
$\P(\tau_n \land t \ge \tau_{n,\dd})$ respectively.

\beg{lem}\label{L3.1} Assume $(A1)$-$(A4)$ and let $K'$ be in $\eqref{D3}$.  Then
$$\P(\tau_n\land t\ge \tau_{n,\dd})\le \ff{\|x-y\|\e^{K't}}{\dd},\ \ t\in [0,T], n\ge 1, x,y\in\H.$$  \end{lem}

\beg{proof} It suffices to prove for $\dd>\|x-y\|.$ By \eqref{D3}, $\e^{-K't}\|X_t^n-Y_t^n\|$ is a supermartingale. So,
$$\dd \e^{-K't} \P(\tau_n\land t\ge \tau_{n,\dd})\le \E\big(\e^{-K'(t\land\tau_n\land\tau_{n,\dd})}\|X^n_{t\land\tau_n\land\tau_{n,\dd}}- Y^n_{t\land\tau_n\land\tau_{n,\dd}}\|\big)\le \|x-y\|.$$ This completes the proof.\end{proof}

Next, we go to estimate $\E(\tau_n\land\tau_{n,\dd}\land t).$  By using $(A1')$ to replace $(A1)$ in the proof of \eqref{D1}, we obtain
$$\d\|X_t^n-Y_t^n\|\le \Big(K\|X_t^n-Y_t^n\|-\ff{\theta \|X_t^n-Y_t^n\|_Q^\kk}{\|X_t^n-Y_t^n\|^{\kk-r}}+I_n(X_t^n-Y_t^n)\Big)\d t +\d M_t^n.$$
Thus, instead of \eqref{D3}, we have
\beq\label{GD} \d\|X_t^n-Y_t^n\|\le \Big(K'\|X_t^n-Y_t^n\|-\ff{\theta \|X_t^n-Y_t^n\|_Q^\kk}{\|X_t^n-Y_t^n\|^{\kk-r}}\Big)\d t +\d M_t^n,\ \ t<T_n\land T.\end{equation}
So, for any $g\in C^2([0,\dd])$ with $g'\ge 0$ and $g''\le 0$,
\beg{equation*}\beg{split} \d g(\|X_t^n-Y_t^n\|) \le  &g'(\|X_t^n-Y_t^n\|) \Big(K'\|X_t^n-Y_t^n\|- \theta\|X_t^n-Y_t^n\|_Q^{\kk} \|X_t^n-Y_t^n\|^{r-\kk}\Big)\d t \\
&+\ff 1 2 g''(\|X_t^n-Y_t^n\|) \d\<M^n\>_t + g'(\|X_t^n-Y_t^n\|)\d M_t^n,\ \ t\le T\land\tau_n\land\tau_{n,\dd}.\end{split}\end{equation*}
Since for $t\le \tau_n\land T$ we have $\|X_t^n-Y_t^n\|\ge \ff 1 n$ such that $h(n\|X_t^n-Y_t^n\|)=1$,   \eqref{D2} implies
\beg{equation}\label{*D}\beg{split} \d\<M^n\>_t &\ge   \ff{4\<QQ_n^{-1}(X_t^n-Y_t^n),X_t^n-Y_t^n\>^2}{\|X_t^n-Y_t^n\|^2\cdot\|Q_n^{-1}(X_t^n-Y_t^n)\|^2}\,\d t\\
&= \ff {4(\|X_t^n-Y_t^n\|^2-\ff 1 n \<Q_n^{-1}(X_t^n-Y_t^n), X_t^n-Y_t^n\>)^2}{\|X_t^n-Y_t^n\|^2\cdot \|Q_n^{-1} (X_t^n-Y_t^n)\|^2}\d t\\
&\ge \ff {4(\ff 1 2 \|X_t^n-Y_t^n\|^4-\ff 1 {n^2}\|Q_n^{-1}(X_t^n-Y_t^n)\|^2\cdot\| X_t^n-Y_t^n\|^2)}{\|X_t^n-Y_t^n\|^2\cdot \|Q_n^{-1} (X_t^n-Y_t^n)\|^2} \d t\\
&\ge  \Big(\ff{2\|X_t^n-Y_t^n\|^2}{\|X_t^n-Y_t^n\|_Q^2} -\ff 4{n^2}\Big)\d t,\ \ t\le \tau_n\land T\end{split}\end{equation} for some constant $c_1>0.$   Since $g''\le 0$, we arrive at
\beq\label{ED} \beg{split} &\d g(\|X_t^n-Y_t^n\|)\\
 &\le   \Big\{g'(\|X_t^n-Y_t^n\|) \Big(K'\|X_t^n-Y_t^n\|-\theta\|X_t^n-Y_t^n\|_Q^{\kk} \|X_t^n-Y_t^n\|^{r-\kk}\Big)\\
&\qquad  +g''(\|X_t^n-Y_t^n\|)\Big(\ff{\|X_t^n-Y_t^n\|^{2}}{\|X_t^n-Y_t\|_Q^{2}}-\ff 2 {n^2}\Big)\Big\}\d t \\
&\qquad+  g'(\|X_t^n-Y_t^n\|)\d M_t^n,\ \ \qquad t\le T\land\tau_n\land\tau_{n,\dd}.\end{split}\end{equation}
By applying  \eqref{ED} with proper choices of $g$,  we will prove   assertions (1)-(3) of Theorem \ref{T1.1} respectively as follows.

\beg{proof}[Proof of Theorem \ref{T1.1}] Without loss of generality, we   only need to prove for $s=0$. Moreover, by the Markov property, it suffices to prove for $t\in [0,T\land 1]$. Below we prove  assertions (1)-(3) respectively for $s=0$ and $t\in (0, T\land 1].$

(1) Let   $\kk>2r$.  We have
$\vv:= \ff{\kk-2r}{\kk}\in (0,1).$ For any $\dd>0$, take
$$g(s)= s-\ff{s^{1+\vv}}{4\dd^\vv},\ \ s\in [0,\dd].$$ Then
\beq\label{GG1} 1\ge g'(s)=1-\ff{(1+\vv)s^\vv}{4\dd^\vv} \ge \ff 1 2,\ \  g''(s)=-\ff{\vv(1+\vv)}{4\dd^\vv s^{1-\vv}}< 0,\ \ s\in (0,\dd].\end{equation}
So,  for any $v\in\V$ with $\|v\|\in [\ff 1 n,\dd]$,
\beg{equation}\label{GG2}\beg{split} G_n(v)&:=  g'(\|v\|) \Big(K'\|v\|-\theta\|v\|_Q^{\kk} \|v\|^{r-\kk}\Big)
  +g''(\|v\|)\Big(\ff{\|v\|^{2}}{\|v\|_Q^{2}}-\ff 2 {n^2}\Big)\\
  &\le K'g(\|v\|) -\ff \theta 2 \|v\|_Q^{\kk} \|v\|^{r-\kk}-\ff{\vv(\vv+1)\|v\|^{1+\vv}}{4\dd^\vv\|v\|_Q^{2}}
  +\ff{\vv(\vv+1)}{2\dd^\vv n^{1+\vv}}.\end{split} \end{equation}
  Noting that $\ff{\|v\|_Q^\kk}{\|v\|^{\kk-r}}= \big(\ff{\|v\|_Q^2}{\|v\|^{\vv+1}}\big)^{\ff\kk 2}$, we obtain
    \beg{equation*}\beg{split} G_n(v)&\le K'g(\|v\|) +\ff{1}{\dd^\vv n^{1+\vv}}-  c_1\bigg\{\Big(\ff{\|v\|_Q^2}{\|v\|^{\vv+1}}\Big)^{\ff\kk 2} +\ff{\|v\|^{1+\vv}}{\dd^\vv \|v\|_Q^2}\bigg\}\\
  & \le K'g(\|v\|)  +\ff{1}{\dd^\vv n^{1+\vv}}- c_1 \inf_{a>0} \big\{a^{-\ff\kk 2} + a \dd^{-\vv}\big\}\\
  &= K'g(\|v\|)  +\ff{1}{\dd^\vv n^{1+\vv}}- c_2\dd^{-\ff{\vv\kk}{\kk+2}},\ \ \ \|v\|\in [n^{-1},\dd]\end{split}\end{equation*} for some constants $c_1,c_2>0.$ Combining this with \eqref{ED}, we arrive at
$$\E(t\land\tau_n\land\tau_{n,\dd}) \le c_3 \dd^{\ff{\vv\kk}{\kk+2}} g(\|x-y\|)+ \ff{c_3}{\dd^{\ff{2\vv}{\kk+2}} n^{1+\vv}},\ \ t\in [0,T\land 1]$$ for some constant $c_3>0.$
This, together with Lemma \ref{L3.1}, yields
\beg{equation}\label{GG3}\beg{split} &\limsup_{n\to\infty} \P(\tau_n>t) \le \limsup_{n\to\infty} \big\{\P(\tau_n\land\tau_{n,\dd}\land t\ge t) +\P(\tau_n\land t\ge \tau_{n,\dd})\big\}\\
&\le \ff{c_3 \dd^{\ff{\vv\kk}{\kk+2}}\|x-y\|}{t} + \ff{\|x-y\|\e^{K'}}\dd,\ \   t\in [0,T\land 1].\end{split}\end{equation}
Therefore, it follows from  Proposition \ref{P2.1}(3)  that
$$\|\nn P_t f\|_\infty \le \osc(f)\inf_{\dd>0} \Big(\ff{c_3 \dd^{\ff{\vv\kk}{\kk+2}}}{t} + \ff{\e^{K'}}\dd\Big)\le
\ff{C\|f\|_\infty}{t^{\ff{\kk+2}{2(\kk+1-r)}}} $$ holds for some constant $C>0$ and all $t\in (0,T].$

(2)  Let $\kk=2r.$   Take $g(s)= \int_0^s \big\{\log(\e+z^{-1})\big\}^{\ff r {1+r}}\d z,\ s\ge 0.$ We have
$$g'(s)=\{\log (\e+s^{-1})\}^{\ff r {1+r}}>0,\ \ \ g''(s)= -\ff{r \{\log (\e +s^{-1})\}^{-\ff 1 {1+r}}} {(1+r)(s+\e s^2)}<0,\ \ s>0.$$ Let $G_n$ be  in \eqref{GG2} and simply take $\dd=1$. Then
\beg{equation*}\beg{split} G_n(v) &\le \Big(K' \|v\|-\ff{\theta \|v\|_Q^{2r}}{\|v\|^r}\Big)\big\{\log(\e +\|v\|^{-1})\big\}^{\ff r {1+r}}
-\ff {r\{\log(\e+\|v\|^{-1})\}^{-\ff 1{1+r}}}{(1+r)(\|v\|+\e\|v\|^2)}\Big(\ff{\|v\|^{2}}{\|v\|_Q^{2}}-\ff 2 {n^2}\Big)\\
&\le K'g(\|v\|) -\ff{\theta\|v\|_Q^{2r}\{\log(\e +\|v\|^{-1})\}^{\ff r {1+r}}}{\|v\|^r}-\ff{r\|v\|\{\log(\e+\|v\|^{-1})\}^{-\ff 1{1+r}}}{(1+r)(1+\e)\|v\|_Q^{2}} +\ff 2 n\\
&\le K'g(\|v\|)+\ff 2 n - \inf_{a>0} \Big(\theta a^{-r} +\ff{r a}{(r+1)(1+\e)}\Big)\\
&= K'g(\|v\|)+\ff 2 n -c_1,\  \ \|v\|\in [n^{-1},1]\end{split}\end{equation*} for some constant $c_1>0.$  Combining this with \eqref{ED} we arrive at
$$\limsup_{n\to\infty} \E(\tau_n\land \tau_{n,1}\land t)\le \ff{g(\|x-y\|)}{c_2},\ \ \ \|x-y\|\in (0, 1), t\le T\land 1$$for some constant $c_2>0$.
Therefore, the first inequality in \eqref{GG3} and Lemma \ref{L3.1} yield
$$\limsup_{n\to\infty} \P(\tau_n>t)\le \ff{c_3}t g(\|x-y\|),\ \ \ \|x-y\|\in (0,1], t\in (0,T\land 1]$$ for some constant $c_3>0.$
Then the proof of (2) is completed by Proposition \ref{P2.1}(3).

(3) Let $\kk\in (r-1, 2r).$ We have $\vv:=\ff{2(\kk+1-r)}{\kk+2}\in (0, 1).$ Take $g(s)= s^\vv,\ s\ge 0.$ Then
 \eqref{GG2} implies
$$ G_n(v) \le \vv K' \|v\|^\vv -  \ff{\vv\theta \|v\|_Q^\kk}{\|v\|^{\kk-r+1-\vv}} -\ff{2\vv(1-\vv)\|v\|^{\vv}}{\|v\|_Q^{2}} +\ff{2\vv(1-\vv)}{n^\vv},\ \ \|v\|\in[n^{-1},1].$$
Since
$$\ff{\|v\|_Q^{\kk}}{\|v\|^{\kk-r+1-\vv}}=\Big(\ff{\|v\|_Q^{2}}{\|v\|^{\vv}}\Big)^{\ff\kk 2},$$ we obtain
$$G_n(v)\le  K'g(\|v\|) +\ff 1 {n^\vv} -\inf_{a>0}\big(\vv\theta a^{-\ff\kk 2} + 2\vv(1-\vv) a\big) =K'g(\|v\|) +\ff 1 {n^\vv}-c_1,\ \ \ \|v\|\in[n^{-1},1]$$ for some constants $c_1>0$.
So,   \eqref{ED} yields
$$\limsup_{n\to\infty} \E(\tau_n\land \tau_{n,1}\land t)\le \ff{g(\|x-y\|)}{c_2}=\ff{\|x-y\|^\vv}{c_2},\ \ \ \|x-y\|\in (0, 1]$$for some constant $c_2>0.$
Combining  this and Lemma   \ref{L3.1} with the first inequality in \eqref{GG3},  we obtain
$$\limsup_{n\to\infty} \P(\tau_n>t)\le \ff{c_3}t g(\|x-y\|),\ \ \ x,y\in\H, t\in (0, T\land 1]$$ for some constant $c_3>0.$
Then the proof is finished by Proposition \ref{P2.1}(3).
\end{proof}

\beg{proof}[Proof of Theorem \ref{T1.2}] (1) Let $r>1.$ Since $(A1')$ is weaker for smaller $\kk$, we assume   $\kk\in (r-1, 2(r-1))$ such that $\vv: = \ff{2(\kk+1-r)}{\kk} \in (0,1).$ Obviously,
\beq\label{VV}  \vv= \ff{2(1+r-\vv)}{\kk+2}-2(1-\vv).\end{equation}
Take
$$g(s)= 1-\e^{-\ll s^\vv} +\gg s^\vv,\ \ s\ge 0,$$ where
  $\ll>1, \gg>0$ will be determined latter on.  Noting that for any $s>0,$
$$g'(s)= \ll \vv s^{\vv-1}\e^{-\ll s^\vv} +\gg \vv s^{\vv-1}> 0,\ \ g''(s) \le -\ll^2\vv^2 s^{2(\vv-1)} \e^{-\ll s^\vv} <0,$$ and $\theta\|\cdot\|_Q^{\kk}\ge \theta'\|\cdot\|^{\kk}$ holds for some constant $\theta'>0$,
we have, for $\|v\|\ge \ff 1 n$,
\beg{equation*}\beg{split} &\ff{\theta g'(\|v\|)\|v\|_Q^\kk}{\|v\|^{\kk-r}}- \Big(\ff{\|v\|^2}{\|v\|_Q^2}-\ff 2{n^2}\Big)^+g''(\|v\|)\\
&\ge \ll\vv\theta \|v\|^{r+\vv-1} \e^{-\ll\|v\|^\vv}\Big(\ff{\|v\|^{2}_Q}{\|v\|^{2}}\Big)^{\ff\kk 2}+ \theta'\gg\vv \|v\|^{\vv+r-1}
+\Big(\ff{\|v\|^2}{\|v\|_Q^2}-\ff 2{n^2}\Big) \ll^2\vv^2 \|v\|^{2(\vv-1)}\e^{-\ll\|v\|^\vv} \\
&\ge \theta'\gg\vv \|v\|^{\vv+r-1}- \ff{2\ll^2\vv^2 \|v\|^{\vv}\e^{-\ll\|v\|^\vv}}{n^\vv} + \e^{-\ll\|v\|^\vv} \inf_{a>0}\Big(
\ll\vv\theta \|v\|^{r+\vv-1} a^{\ff\kk 2} + \ff{\ll^2\vv^2}{a \|v\|^{2(1-\vv)}}\Big)\\
&= \theta'\gg\vv \|v\|^{\vv+r-1}- \ff{2\ll^2\vv^2 \|v\|^{\vv}\e^{-\ll\|v\|^\vv}}{n^\vv} + c_2 \ll^{\ff{2(\kk+1)}{\kk+2}}
\|v\|^{\ff{2(1+r-\vv)}{\kk+2}-2(1-\vv)}\e^{-\ll\|v\|^\vv}\\
&= \theta'\gg\vv \|v\|^{\vv+r-1}- \ff{2\ll^2\vv^2 \|v\|^{\vv}\e^{-\ll\|v\|^\vv}}{n^\vv} + c_2 \ll^{\ff{2(\kk+1)}{\kk+2}} \|v\|^{\vv}\e^{-\ll\|v\|^\vv}
\end{split}\end{equation*} for some constant $c_2>0$, where the last step is due to \eqref{VV}.
Therefore,
\beg{equation*}\beg{split} G_n(v) & := g'(\|v\|)\Big(K'\|v\| -\ff{\theta \|\cdot\|_Q^{\kk}}{\|v\|^{\kk-r}}\Big)+
  \Big(\ff{c_1 \|v\|^2}{\|v\|_Q^2}-\ff 2{n^2}\Big) g''(\|v\|)\\
&\le K'\gg \vv \|v\|^\vv- \theta'\gg\vv \|v\|^{\vv+r-1} -\Big(c_2\ll^{\ff{2(\kk+1)}{\kk+2}}
-K'\ll\vv -\ff{2\ll^2\vv^2}{n^\vv}\Big)\|v\|^{\vv}\e^{-\ll\|v\|^\vv}\end{split}\end{equation*} holds for $\|v\|\ge \ff 1 n.$
Obviously, there exist $\ll,n_0\ge 1$ such that
$$c_2\ll^{\ff{2(\kk+1)}{\kk+2}}
-K'\ll\vv -\ff{2\ll^2\vv^2}{n^\vv_0}\ge \ff {c_2} 2\ll^{\ff{2(\kk+1)}{\kk+2}} =:c_3>0.$$ Moreover, take
$$\gg= \ff{c_3\exp[-\ll(\ff{2K'}{\theta'})^{\ff\vv{r-1}}]}{2K'\vv}$$ such that
$$2K'\gg\vv \|v\|^\vv \le \beg{cases} \theta'\gg\vv \|v\|^{r+\vv-1},\ &\text{if} \ \|v\|\ge (\ff{2K'}{\theta'})^{\ff 1{r-1}},\\
c_3 \|v\|^\vv \e^{-\ll \|v\|^\vv}, &\text{if}\ \|v\|< (\ff{2K'}{\theta'})^{\ff 1{r-1}}.\end{cases}$$
Therefore, there exists a constant $c_4>0$ such that
$$G_n(v)\le -\gg K'\vv \|v\|^\vv\le - c_4 g(\|v\|), \ \ n\ge n_0, \|v\|\ge\ff 1 n.$$ Combining this with \eqref{ED} we obtain
$$\e^{c_4 t}\E\big(g(\|X_t^n-Y_t^n\|)1_{\{t\le\tau_n\}}\big)\le \E\big\{\e^{c_4(t\land\tau_n)}g(\|X_{t\land\tau_n}^n-Y_{t\land\tau_n}^n\|)\big\}\le g(\|x-y\|).$$ So, there exists a constant $c_5>0$ such that
$$\E\big(\|X_t^n-Y_t^n\|^\vv 1_{\{t\le\tau_n\}}\big)\le  c_5 \|x-y\|^\vv \e^{-c_4t}.$$
On the other hand,
by Theorem \ref{T1.1},   there exists $p\ge 1$ such that
$$|P_{t, t+1} f(x)-P_{t,t+1}f(y)|^p \le c \|f\|_\infty^p (\|x-y\|^\vv\land 1),\ \ x,y\in \H, t\ge 0, f\in \B_b(\H)$$ holds for some constant $c>0.$ Thus, for $\|f\|_\infty\le 1$,
\beg{equation*}\beg{split} &|P_{t+1}f(x)-P_{t+1}f(y)|^p = |\E \{P_{t,t+1}f(X_t^n)-P_{t,t+1}f(Y_t^n)\}|^p\\
 &\le \E |P_{t,t+1}f(X_t^n)-P_{t,t+1}f(Y_t^n)|^p
 \le c  \E (1\land \|X_t^n-Y_t^n\|^\vv)\\
 & \le cc_5 \|x-y\|^\vv \e^{-c_4t} +c \big\{\E(\|X_t^n-Y_t^n\|1_{\{t>\tau_n\}})\big\}^{\vv}\\
&\le cc_5 \|x-y\|^\vv \e^{-c_4t}+ c\Big(\ff{\e^{K't}}n\Big)^{\vv},\end{split}\end{equation*} where the last step is due to \eqref{*W}.
Letting $n\to\infty$, we get
$$|P_{t+1}f(x)-P_{t+1}f(y)|^p\le cc_5 \|x-y\|^\vv \e^{-c_4t},\ \ t\ge 0, x,y\in\H.$$ Considering the equation \eqref{E1} from time $s$ rather than $0$, this inequality becomes
$$|P_{s,s+t+1}f(x)-P_{s,s+t+1}f(y)|^p\le cc_5 \|x-y\|^\vv \e^{-c_4t},\ \ t\ge 0, x,y\in\H.$$ Thus,
as explained above, we have
\beq\label{PP}\beg{split} |P_{t+2}f(x)-P_{t+2}f(y)|^p&= \big|\E\big[(P_{1,t+2}f)(X_1(x))-(P_{1,t+2}f)(X_1(y))\big]\big|^p\\
 &\le cc_5 \e^{-c_4t}\E \|X_{1}(x)-X_{1}(y)\|^\vv ,\  \ t\ge 0, x,y\in\H.\end{split}\end{equation} By $(A1')$ and It\^o's formula,
$$\d \|X_{s}(x)- X_{s}(y)\|^2\le \big\{c_6- k  \|X_{s}(x)- X_{s}(y)\|^{1+r}\big\}\d s +\d M_s, s\ge 0$$ holds for some constants $c_6,k>0$ and some martingale
$M_s$. Since $r>1$, this implies
$$  \sup_{x,y\in\H}\E \|X_{1}(x)- X_{1}(y)\|^2<\infty.$$ Therefore, \eqref{PP} implies
$$|P_{t+2}f(x)-P_{t+2}f(y)|^p \le C \e^{- c_4 t},\ \ t\ge 0, \|f\|_\infty\le 1$$ for some constant $C>0$,  so that the proof of (1) is finished.

(2) Let $r=1.$ It suffices to prove for large $t>0.$ By It\^o's formula and \eqref{EC}, we have
$$\E\|X_t(x)-X_t(y)\|^2\le \|x-y\|^2\e^{-2Kt}.$$ Therefore, for any $p\in (0,2]$,
$$\E \|X_t(x)-X_t(y)\|^p \le \|x-y\|^p\e^{-pKt},$$ and, since $0\le s\mapsto \int_0^s \ss{\log(\e +   z^{-1})}\d z$ is concave,
\beg{equation*}\beg{split} \E\int_0^{\|X_t(x)-X_t(y)\|}\ss{\log\big(\e+ z^{-1}\big)}\d z&\le \int_0^{\|x-y\|\e^{-Kt}}\ss{\log\big(\e+ z^{-1}\big)}\d z\\
&\le c \|x-y\|\e^{-Kt}\ss{\log\big(\e+ \|x-y\|^{-1}\e^{Kt}\big)}\end{split}\end{equation*} for some constant $c>0.$  Combining these with Theorem \ref{T1.1} and noting that
$$|P_{t+1}f(x)-P_{t+1}f(y)|\le \E |P_{t,t+1}f(X_t(x))- P_{t,t+1}f(X_t(y))|,$$ we prove (2).
\end{proof}

\section{Proof of Theorem \ref{T1.3}}

According to condition $(A1'')$, we will need to estimate moments on $\|X_t\|_\V$. To this end, we first introduce the following lemma  which is implied by  $(A3)$ and  It\^o's formula for $\|X_s\|^2,$ see \cite{L,Wbook} for $B=0$.  To save space, we omit the proof.

\beg{lem}\label{L5.1} Assume $(A1)$-$(A4)$. Then for any $p>0$  there exists a constant $c(p)>0$ such that
$$\E\bigg(\int_0^t \|X_s(x)\|_\V^{1+r}\d s\bigg)^p \le c(p)(1+\|x\|^2)^p,\ \ x\in\H, t\in [0,T\land 1].$$
If moreover $B$ is bounded, then there exist two constants $\ll,c>0$ such that
$$\E\e^{\ll\int_0^t\|X_s(x)\|_{\V}^{1+r}\d s}\le \e^{c(1+\|x\|^2)},\ \ x\in\H, t\in [0,T\land 1].$$
\end{lem}

\

Similarly to the proof of \eqref{GD} using $(A1')$,  it is easy to see that $(A1'')$ implies
$$\d \|X_t^n-Y_t^n\|\le \bigg(K' \|X_t^n-Y_t^n\|-\ff{\theta\|X_t^n-Y_t^n\|_Q^\kk}{\|X_t^n-Y_t^n\|^{\kk-1} h_t^{1-r}}\bigg)\d t + \d M_t^n,\ \ t<T_n\land T$$ for some constant $K'>0$, $M_t^n$ in \eqref{D2},  and
$$h_t:= \|X_t^n\|_\V\lor \|Y_t^n\|_\V.$$ So, for any $\dd>0$ and $g\in C^2([0,\dd])$ with $g'\ge 0$ and $g''\le 0$, we have
\beq\label{C*0} \beg{split}&\d g(\|X_t^n-Y_t^n\|) \le H_t(X_t^n-Y_t^n)\d t + g'(\|X_t^n-Y_t^n\|) \d M_t^n,\ \ t<T_n\land T,\\
&H_t(v):= g'(\|v\|) \bigg(K' \|v\|-\ff{\theta\|v\|_Q^\kk}{\|v\|^{\kk-1} h_t^{1-r}}\bigg)+\ff{g''(\|v\|)} 2 \bigg(\ff{\|v\|^2}{\|v\|_Q^2}-\ff 2 {n^2}\bigg),\ \  0\ne v\in\V.\end{split} \end{equation}

\beg{proof}[Proof of Theorem \ref{T1.3}] As explained in the proof of Theorem \ref{T1.1} that we only prove this theorem for $s=0$ and $t\le T\land 1$. Below we prove  assertions (1)-(3) respectively by \eqref{C*0} with different choices of $g$.

(a) Let $\kk>2$. We have $\vv:=\ff{\kk-2}\kk\in (0,1).$ For any $\dd>0$, take
$$g(s)= s-\ff{s^{1+\vv}}{4\dd^\vv},\ \ s\in [0,\dd].$$ Obviously,
\beq\label{C*2} \ff s 2 \le g(s)\le s,\ \ \ff 1 2 \le g'(s)\le 1,\ \ g''(s)= -\ff{\vv(1+\vv)s^{\vv-1}}{4\dd^\vv},\ \ s\in (0, \dd].\end{equation}
Then letting
$$A_s =\ff{\|X_s^n-Y_s^n\|^{1+\vv}}{\|X_s^n-Y_s^n\|_Q^2},\ \ \ B_s=\ff{\|X_s^n-Y_s^n\|_Q^\kk}{\|X_s^n-Y_s^n\|^{\kk-1}},$$ we obtain
$$H_t(X_s^n-Y_s^n) \le 2 K' g(\|X_s^n-Y_s^n\|) +\ff 1{\dd^\vv n^{1+\vv}} -c_1 \big(B_s h_s^{r-1} + \dd^{-\vv} A_s\big),\ \ s\le T\land \tau_n\land\tau_{n,\dd}$$ for some constant $c_1>0.$ Thus,   it follows from \eqref{C*0} that
\beq\label{C*3} \E \int_0^{t\land\tau_n\land\tau_{n,\dd}} \big(B_s h_s^{r-1} + \dd^{-\vv} A_s\big)\d s
\le c_2 \Big(\|x-y\|+\ff 1{\dd^\vv n^{1+\vv}}\Big),\ \ t\le T\land 1,\dd>0 \end{equation} holds for some constant $c_2>0$.
It is easy to see that $B_s= A_s^{-\ff\kk 2}.$ So,
$$B_s h_s^{r-1} + \dd^{-\vv}A_s\ge \inf_{a>0}\{a^{-\ff\kk 2} h_s^{r-1} + \dd^{-\vv} a\} =c_3 \dd^{\ff{2-\kk}{2+\kk}} h_s^{-\ff{2(1-r)}{2+\kk}}$$ holds for some constant $c_3>0.$
Combining this with \eqref{C*3}, we obtain
\beq\label{C*4} \limsup_{n\to\infty} \E\int_0^{t\land \tau_n\land\tau_{n,\dd}}h_s^{-\ff{2(1-r)}{2+\kk}}\d s\le c_4\|x-y\|\dd^{\ff{\kk-2}{2+\kk}}   \end{equation}for some constant
$c_4>0.$  Then for any $R>0$,
\beq\label{CN} \beg{split} &\limsup_{n\to\infty}\P(\tau_n\land\tau_{n,\dd}\ge t)
  \le \limsup_{n\to\infty}\P\bigg(\int_0^{t\land \tau_n\land\tau_{n,\dd}}h_s^{-\ff{2(1-r)}{2+\kk}}\d s\ge \int_0^{t}h_s^{-\ff{2(1-r)}{2+\kk}}\d s\bigg)\\
&\le \limsup_{n\to\infty}\P\bigg(\int_0^{t\land \tau_n\land\tau_{n,\dd}}h_s^{-\ff{2(1-r)}{2+\kk}}\d s\ge\ff{ t^{\ff{4+\kk+\kk r}{(2+\kk)(r+1)}}}{(\int_0^{t}h_s^{1+r}\d s)^{\ff{2(1-r)}{(2+\kk)(1+r)}}}\bigg)\\
&\le \limsup_{n\to\infty}\bigg\{\P\bigg(\int_0^{t\land \tau_n\land\tau_{n,\dd}}h_s^{-\ff{2(1-r)}{2+\kk}}\d s\ge\ff{ t^{\ff{4+\kk+\kk r}{(2+\kk)(r+1)}}} R\bigg)
 + \P\bigg(\int_0^{t}h_s^{1+r}\d s> R^{\ff{(2+\kk)(1+r)}{2(1-r)}}\bigg)\bigg\}\\
&\le \ff{c_4\|x-y\|\dd^{\ff{\kk-2}{\kk+2}}R}{t^{\ff{4+\kk+\kk r}{(2+\kk)(r+1)}}} + \limsup_{n\to\infty}\P\bigg(\int_0^{t}h_s^{1+r}\d s> R^{\ff{(2+\kk)(1+r)}{2(1-r)}}\bigg).
\end{split}\end{equation} By Proposition \ref{P2.1}(1) and Lemma \ref{L5.1}, for any $p>0$ there exists a constant $c(p)>0$ such that
\beq\label{CN2}  \P\bigg(\int_0^{t}h_s^{1+r}\d s> R^{\ff{(2+\kk)(1+r)}{2(1-r)}}\bigg)\le \ff{c(p)(1+\|x\|^2+\|y\|^2)^{\ff{2p(1-r)}{(2+\kk)(r+1)}}}{R^p},\end{equation} and when $B$ is bounded there exists a constant $c>0$ such that
\beq\label{CN3} \P\bigg(\int_0^{t}h_s^{1+r}\d s> R^{\ff{(2+\kk)(1+r)}{2(1-r)}}\bigg)\le\exp\Big[c(1+\|x\|^2+\|y\|^2) -\ll R^{\ff{(2+\kk)(1+r)}{2(1-r)}}\Big].\end{equation}
Now, for any $p>0$ it follows from \eqref{CN} and \eqref{CN2} that
\beg{equation*} \beg{split} \limsup_{n\to\infty}\P(\tau_n\land\tau_{n,\dd}\ge t) &\le \inf_{R>0}
\bigg\{\ff{c_4\|x-y\|\dd^{\ff{\kk-2}{\kk+2}}R}{t^{\ff{4+\kk+\kk r}{(2+\kk)(r+1)}}} + \ff{c(p)(1+\|x\|^2+\|y\|^2)^{\ff{2p(1-r)}{(2+\kk)(r+1)}}}{R^p}\bigg\}\\
&=\ff{c_5 (1+\|x\|^2+\|y\|^2)^{\ff{2p(1-r)}{(p+1)(2+\kk)(r+1)}}\|x-y\|^{\ff p{p+1}} \dd^{\ff{p(\kk-2)}{(p+1)(\kk+2)}}}{t^{\ff{p(4+\kk+\kk r)}{(p+1)(2+\kk)(r+1)}}}\end{split}\end{equation*}
for some constant $c_5>0.$ Combining this with Lemma \ref{L3.1}, we obtain
\beg{equation*}\beg{split} &\limsup_{n\to\infty} \P(\tau_n\ge t)  \le \limsup_{n\to\infty} \big\{\P(\tau_n\land\tau_{n,\dd}\ge t) +\P(\tau_n\land t\ge \tau_{n,\dd})\big\}\\
&\le c_6\bigg( \ff{\|x-y\|}\dd +\ff{(1+\|x\|^2+\|y\|^2)^{\ff{2p(1-r)}{(p+1)(2+\kk)(r+1)}}\|x-y\|^{\ff p{p+1}} \dd^{\ff{p(\kk-2)}{(p+1)(\kk+2)}}}{t^{\ff{p(4+\kk+\kk r)}{(p+1)(2+\kk)(r+1)}}}\bigg),\  \ \dd>0\end{split}\end{equation*} for some constant $c_6>0.$ Minimizing the upper bound in $\dd>0$, we arrive at
$$\limsup_{n\to\infty} \P(\tau_n\ge t) \le c_7 (1+\|x\|^2+\|y\|^2)^{\ff{2p(1-r)}{(2p\kk+\kk+2)(r+1)}}\|x-y\|^{\ff {2p\kk}{2p\kk+\kk+2}} t^{-\ff{p(4+\kk+\kk r)}{(2p\kk+\kk+2)(r+1)}}$$
for some constant $c_7>0$. According to Proposition \ref{P2.1}(1), this implies the first assertion in (1).

When $B$ is bounded, by \eqref{CN} and \eqref{CN3} we have
$$\limsup_{n\to\infty}\P(\tau_n\land\tau_{n,\dd}\ge t)\le
 \ff{c_4\|x-y\|\dd^{\ff{\kk-2}{\kk+2}}R}{t^{\ff{4+\kk+\kk r}{(2+\kk)(r+1)}}} + \exp\Big[c(1+\|x\|^2+\|y\|^2)-\ll R^{\ff{(1+r)(2+\kk)}{2(1-r)}}\Big] $$ for any $R>0$.
 Taking
 $$R= \Big\{\ff 1 \ll \big(c(1+\|x\|^2+\|y\|^2)+\log(\e+\|x-y\|^{-1})\big)\Big\}^{\ff{2(1-r)}{(1+r)(2+\kk)}},$$
 we arrive at
 $$\limsup_{n\to\infty}\P(\tau_n\land\tau_{n,\dd}\ge t)\le c_8 \dd^{\ff{\kk-2}{\kk+2}} \|x-y\| \big(1+\|x\|^2+\|y\|^2+\log(\e+\|x-y\|^{-1})\big)^{\ff{2(1-r)}{(1+r)(2+\kk)}}$$
 for some constant $c_8>0.$  Combining this with Lemma \ref{L3.1} and Proposition \ref{P2.1}(3), we prove the second assertion in (1).

(b) Let $\kk=2.$ Take
$$g(s)= \int_0^s \ss{\log (\e + z^{-1})}\,\d z,\ \ s\ge 0.$$
By \eqref{C*0} we have
$$H_t(v)\le K'g(\|v\|) +\ff 1 n -\ff{\theta \|v\|_Q^2\ss{\log (\e +\|v\|^{-1})}}{\|v\| h_t^{1-r}}
-\ff{\|v\|}{4(\e +1) \|v\|_Q^2\ss{\log (\e +\|v\|^{-1})}}$$ for $\|v\|\in [n^{-1},1].$ So, according to \eqref{C*0},
\beq\label{C*4} \E \int_0^{t\land\tau_n\land\tau_{n,1}} \big(A_s^{-1} h_s^{r-1} + A_s\big)\d s\le c_1 \Big(g(\|x-y\|)+ \ff 1 n\Big),\ \ t\in [0,T\land 1]\end{equation} holds for some constant $c_1>0$ and
$$A_s:= \ff{\|X_s^n-Y_s^n\|}{\|X_s^n-Y_s^n\|_Q^2 \ss{\log (\e +\|X_s^n-Y_s^n\|^{-1})}}.$$
Since $A_s^{-1}h_s^{r-1}+A_s\ge 2 h_s^{\ff{r-1}2},$ this implies
$$\E\int_0^{t\land\tau_n\land\tau_{n,1}} h_s^{\ff{r-1}2}\d s \le \ff{c_1} 2 \big\{g(\|x-y\|)+n^{-1}\big\}.$$
Similarly to \eqref{CN} with $\dd=1$ and $\kk=2$, this implies
$$\limsup_{n\to\infty}\P(\tau_n\land\tau_{n,1}\ge t)\le\ff{c_2g(\|x-y\|)R}{t^{\ff{3+r}{1-r}}} + \P\bigg(\int_0^t h_s^{1+r}\d s\ge R^{\ff{2(1+r)}{1-r}}\bigg),\ \ R>0$$for some constant $c_2>0.$ Combining this with \eqref{CN2} and \eqref{CN3}, we conclude that for any $p>0$ there exists a constant $c(p)>0$ such that
\beq\label{CNa} \limsup_{n\to\infty}\P(\tau_n\land\tau_{n,1}\ge t)\le \ff{c(p)(1+\|x\|^2+\|y\|^2)^{\ff{p(1-r)}{2(1+p)(1+r)}}g(\|x-y\|)^{\ff p{p+1}}}{t^{\ff{p(3+r)}{2(1+p)(1+r)}}}\end{equation} holds, and when $B$ is bounded
\beq\label{CNb} \limsup_{n\to\infty}\P(\tau_n\land\tau_{n,1}\ge t)\le  c g(\|x-y\|)\big(1+\|x\|^2+\|y\|^2+\log(\e+\|x-y\|^{-1})\big)^{\ff{1-r}{2(1+r)}}\end{equation}holds for some constant
$c>0.$ Therefore, the assertions in  (2) follow    from Lemma \ref{L3.1} with $\dd=1$ and Proposition \ref{P2.1}(3).

(c) Let $\kk\in (0,2).$  We have $\vv:= \ff{2\kk}{\kk+2}\in (0,1).$ Take $g(s)= s^\vv$. By \eqref{C*0} we have
$$H_t(v)\le c_1 \Big(g(\|v\|)+\ff 1 {n^\vv}\Big) - c_2 \Big(\ff{\|v\|_Q^\kk}{\|v\|^{\kk-\vv}h_t^{1-r}} + \ff{\|v\|^\vv}{\|v\|_Q^2}\Big),\ \ v\ne 0$$ for some constants $c_1,c_2>0$. Let
$$A_s =\ff{\|X_s^n-Y_s^n\|^\vv}{\|X_s^n-Y_s^n\|_Q^2}.$$ By the choice of $\vv$ we have
$$ \ff{\|X_s^n-Y_s^n\|_Q^\kk}{\|X_s^n-Y_s^n\|^{\kk-\vv}}= A_s^{-\ff\kk 2},$$ so that \eqref{C*0} yields
\beq\label{C*5} \E \int_0^{t\land\tau_n\land\tau_{n,1}} \big(A_s^{-\ff\kk 2} h_s^{r-1} + A_s\big)\d s\le c_3 \Big(g(\|x-y\|)+ \ff 1 {n^\vv}\Big),\ \ t\le T\land 1\end{equation} for some constant $c_3>0.$  As explained in (b), this implies \eqref{CNa},  and also \eqref{CNb} when $B$ is bounded, for   $g(\|x-y\|)= \|x-y\|^{\ff{2\kk}{\kk+2}}$. Therefore, the assertions in (3)
from Lemma \ref{L3.1} with $\dd=1$ and Proposition \ref{P2.1}(3).\end{proof}

\section{Applications to specific models  } In this section we apply  Theorems \ref{T1.1} and \ref{T1.2} to the stochastic generalized porous media equations and the stochastic $p$-Laplace equations, and apply Theorem \ref{T1.3} as well as Proposition \ref{P*} to the stochastic generalized fast-diffusion equations.

\subsection{Stochastic generalized porous media equations}

Let $(E,\B,{\m})$ be a separable probability space and
$(L,\D(L))$ a negative definite self-adjoint linear operator on
$L^2({\bf m})$ having
 discrete spectrum.  Let $$
  (0<) \ll_1\le \ll_2\le \cdots
$$ be all eigenvalues of $-L$ including multiplicities with unit eigenfunctions $\{e_i\}_{i\ge 1}$.
Let $\H$ be the dual space of the
$\D((-L)^{\ff 12})$ with respect to $L^2(\m)$; i.e. $\H$ is the  completion of $L^2(\m)$ under the inner
product
$$\<x,y\>:= \sum_{i=1}^\infty \ff 1 {\ll_i}
\m(xe_i)\m(ye_i),$$ where $\m(x):= \int_E x\d\m$ for $x\in L^1(\m).$    Let
$$\Psi, \Phi: [0,\infty)\times \R \to \R$$
be   measurable,  and be continuous in the second variable.
   We consider the equation
\beq\label{2.1.1} \d X_t = \big\{L\Psi(t,X_t)+ \Phi(t,X_t)\big\}\d t
+B(t,X_t)\d W^{(1)}_t + Q\d W_t^{(2)},\end{equation} where $W_t^{(1)}$ and $W_t^{(2)}$ are independent cylindrical Brownian motions on $\H$, $Q\in \L_{HS}(\H)$ and $B: [0,\infty)\times \H  \to \L_{HS}(\H)$ is measurable.

To verify conditions $(A1')$ and $(A2)$-$(A4)$ for $A(t,v):= L \Psi(t,X_t)+ \Phi(t,X_t),$ we assume
that for a fixed constant $r\ge 1$,
\beq\label{4.1}\beg{split} & |\Psi(t,s)|+|\Phi(t,s)|\le c(1+|s|^r), \ \ s\in\R, t\ge 0,\\
&\ff 1 2 \|B(t,x)-B(t,y)\|_{HS}^2 -\m\big((\Psi(t,x)-\Psi(t,y))( x-y)\big) \\
&\quad +\m\big((\Phi(t,x)-\Phi(t,y)) (-L)^{-1} (x-y)\big)\le
K\|x-y\|^2-\theta \|x-y\|_{r+1}^{r+1},\ \  t\ge 0\end{split} \end{equation} holds for some constants $c, K,\theta>0$ and all $x,y\in L^{r+1}(\m),$
where $\|\cdot\|_{r+1}$ is the norm in $L^{1+r}(\m).$  Obviously, this condition is satisfied provided
\beq\label{LIP}\|B(t,x)-B(t,y)\|_{HS}\le c_0 \|x-y\|,\ \ x,y\in\H \end{equation} holds for some constant $c_0>0,$
$\Psi(t,s)= h(t) s^r$ and $\Phi(t,s)= g(t) s$ with
$0<\inf h\le \sup h<\infty$ and $ \|g\|_\infty<\infty$, where $s^r:= |s|^{r-1}s$ for $s\in\R.$

\

Now, let $\V=L^{1+r}(\m)$ and let $\V^*$ be the dual space of $\V$ with respect to $\H$. Then it is easy to see that \eqref{4.1} implies
$(A2)$-$(A4)$ for
$$A(t,v):= L\Psi(t,v) +\Phi(t,v).$$   Therefore, to apply Theorems \ref{T1.1} and \ref{T1.2}, it remains to verify $(A1')$.

 \beg{prp}\label{P4.1}  Assume $\eqref{4.1}$. Let $Qe_i= q_i e_i, i\ge 1$, where $\{q_i\}_{i\ge 1}\subset \R$ satisfy $\sum_{i=1}^\infty q_i^2<\infty.$ If
 for some $\kk\ge 1+r$ \beq\label{*E} \sup_{i\ge 1}  \ll_i^{-1}q_i^{-\ff {2\kk}{1+r}}<\infty,\end{equation}
 then assertions in Theorems $\ref{T1.1}$ and $\ref{T1.2}$ hold.
 \end{prp}
 \beg{proof} Note that   $\sum_{i=1}^\infty q_i^2<\infty$ ensures $Q\in \L_{HS}(\H)$. By the second inequality in \eqref{4.1},
  \beq\label{M1}  \beg{split} &  _{\V^*}\<A(t, v_1)-A(t, v_2), v_1-v_2\>_\V + \ff 1 2  \|B(t, v_1)-B(t,v_2)\|_{HS}^2\\
   &\le K_1|v_1-v_2\|^2-\theta_1\|v_1-v_2\|^{r+1}_\V,\  \ t\ge 0,  v_1,v_2\in  \V:=L^{1+r}(\m) \end{split}\end{equation} holds for some constants $K_1,\theta_1>0.$
On the other hand,    \eqref{*E} implies
 \beg{equation*}\beg{split} \|x\|_Q^2& :=\sum_{i\ge 1} q_i^{-2} \ll_i^{-1}\m(e_i x)^2\le \bigg(\sum_{i\ge 1} \ll_i^{-1}\m(e_i x)^2\bigg)^{\ff{\kk-r-1}\kk}
 \bigg(\sum_{i\ge 1} q_i^{-\ff {2\kk}{1+r}} \ll_i^{-1}\m(e_ix)^2\bigg)^{\ff{1+r}\kk}\\
 &\le C\|x\|^{\ff{2(\kk-1-r)}\kk}\|x\|_2^{\ff{2(1+r)}\kk}\le C \|x\|^{\ff{2(\kk-1-r)}\kk}\|x\|_\V^{\ff{2(1+r)}\kk}  \end{split}\end{equation*}
for some constant $C>0$. Combining this with \eqref{M1}, we prove  $(A1')$.  \end{proof}

 Below, we present a simple example to illustrate this result.

 \paragraph{Example  4.1.}\ Let $r>1$ and $\DD$ be the Dirichlet Laplacian on a bounded domain $D\subset \R^d$. Let $L= -(-\DD)^\gg$ for some constant $\gg>0$. Let $\m$ be the normalized Lebesgue measure on $D$. Take $$\Phi(s)=cs,\ \  \Psi(t,s)=s^r:= s|s|^{r-1}, \ \ Q  e_i= c_ii^{-\dd} e_i,\  i\ge 1$$ for some constants $c\in\R$,   $\dd>\ff 1 2 $ and $\{c_i\}_{i\ge 1}$ with $0<\inf |c_i|\le \sup |c_i|<\infty$. Moreover, let
 If $\gg\ge \dd d$, then assertions in Theorem \ref{T1.1} and Theorem \ref{T1.2}(1) hold for $\kk:= \ff{\gg(1+r)}{\dd d}\ge 1+r.$

 \beg{proof} Obviously, for the specific functions $\Phi$ and $\Psi$,   \eqref{LIP} implies \eqref{4.1}.     Let $q_i:= c_ii^{-\dd}$. It is easy to see that   $\sum_{i\ge 1}q_i^2<\infty$ and \eqref{4.1} holds.  Next, we have
 $$\ll_i\ge ci^{\ff{2\gg}d},\ \ i\ge 1$$ for some constant $c>0$. Then \eqref{*E} holds for $\kk:=  \ff{\gg(1+r)}{\dd d}.$ So, the conclusion follows from  Proposition \ref{P4.1}.\end{proof}

\subsection{Stochastic $p$-Laplacian equations}

Let $D\subset \R^d$ be an open domain, let $\m$ be the normalized volume measure on $D$,  and let $p\ge 2$ be a constant. Let $\H_0^{1,p}(D)$ be the closure of $C_0^\infty(D)$ with respect to the norm
$$\|f\|_{1,p}:= \|f\|_p  +\|\nn f\|_p,$$ where $\|\cdot\|_p$ is the norm in $L^p(\m)$.  Let $\H=L^2(\m)$ and $\V=\H_0^{1,p}(D)$. By the $L^p$-Poincar\'e inequality, there exists a constant $C>0$ such that
$\|f\|_{1,p}\le C \|\nn f\|_{L^p(\m)}.$ Consider the SPDE
$$\d X_t = {\rm div}\big(|\nn X_t|^{p-2}\nn X_t\big) \d t +B(t,X_t)\d W^{(1)}_t + Q\d W_t^{(2)},$$ where $W_t^{(1)}$ and $W_t^{(2)}$ are independent cylindrical Brownian motions on $\H$, $Q\in \L_{HS}(\H)$ and $B: [0,\infty)\times \H  \to \L_{HS}(\H)$ is measurable such that \eqref{LIP} holds for some constant $c_0>0$.  Then $(A2)$-$(A4)$ hold for (see \cite[Example 3.3]{L})
$$A(t,v):= {\rm div} (|\nn v|^{p-2}\nn v).$$
Moreover,  by \cite[Lemma 3.1]{L} and \eqref{LIP},  \eqref{M1} holds for some constants $K_1,\theta_1>0.$

To verify $(A1')$, we simply consider $d=1$ and   $D=(0,1).$   Let $\DD$ be the Dirichlet Laplacian on $(0,1)$, then $\{(\pi i)^2\}_{i\ge 1}$ are all eigenvalues of $-\DD$ with unit eigenfunctions $e_i(x):= \ss 2 \sin (i\pi x)$.

\beg{prp}\label{P4.2} Let $D=(0,1)\subset \R$ and $B$ satisfy $\eqref{LIP}$. Let $Qe_i= q_i e_i, i\ge 1,$ where $\{q_i\}_{i\ge 1}\subset \R$ satisfy $\sum_{i=1}^\infty q_i^2<\infty.$  If there exists $\kk\ge p$ such that
\beq\label{**E} \sup_{i\ge 1} q_i^{-1} i^{-\ff p\kk}<\infty,\end{equation} then assertions in Theorem $\ref{T1.1}$ and Theorem $\ref{T1.2}(1)$ hold for $r:=p-1.$ When $p=2\ ($i.e. $r:= p-1=1)$ and   $K:=\pi^2-\ff 1 2c_0^2>0,$ Theorem $\ref{T1.2}(2)$ holds. \end{prp}

\beg{proof} Let $r= p-1$. Obviously,   \eqref{**E} and Schwartz's inequality imply
\beg{equation*}\beg{split} \|x\|_Q^2 & =\sum_{i\ge 1} q_i^{-2} \m(xe_i)^2 \le \bigg(\sum_{i\ge 1} \m(xe_i)^2\bigg)^{\ff{\kk-p}\kk}\bigg(\sum_{i\ge 1} q_i^{-\ff {2\kk}p}\m(xe_i)^2\bigg)^{\ff p\kk}\\
&\le \|x\|_2^{\ff{2(\kk-p)}\kk}\bigg(\sum_{i\ge 1} (\pi i)^2 \m(xe_i)^2\bigg)^{\ff p\kk} \sup_{i\ge 1} q_i^{-2} (\pi i)^{-\ff{2p}\kk}\\
&\le C \|x\|^{\ff{2(\kk-p)}\kk} \m(|\nn x|^2)^{\ff p\kk}   \le C \|x\|^{\ff{2(\kk-1-r)}\kk} \|x\|_\V^{\ff{2(1+r)}\kk}  \end{split}\end{equation*} for some constant $C>0.$ Then $(A1')$ follows from \eqref{M1}, which is implied by \eqref{LIP} and \cite[Lemma 3.1]{L}.  It remains to show that when $r=1$ (i.e. $p=2$), the condition \eqref{EC} in Theorem \ref{T1.2}(2) holds.
Indeed,  when $p=2$, the integration by parts formula, the Poincar\'e inequality and  \eqref{LIP} yield
\beg{equation*}\beg{split}& _{\V^*}\<A(t,v_1)-A(t,v_2), v_1-v_2\>_\V +\ff 1 2\|B(t,v_1)-B(t,v_2)\|_{HS}^2 \\
&= - \m(|\nn(v_1-v_2)|^2) +\ff 1 2\|B(t,v_1)-B(t,v_2)\|_{HS}^2\\
&\le -\m(|\nn(v_1-v_2)|^2) +\ff {c_0^2}2 \|v_1-v_2\|^2 \\
&\le -\Big(\pi^2-\ff {c_0^2}2 \Big)\|v_1-v_2\|^2,\ \ v_1,v_2\in \V.\end{split}\end{equation*}
\end{proof}

\paragraph{Example  4.2.}\ In the situation of Proposition \ref{P4.2}, take   $q_i=c_ii^{-\dd} $ for some constants   $\dd\in (\ff 1 2, 1]$ and $\{c_i\}_{i\ge 1}$ with $0<\inf |c_i|\le \sup |c_i|<\infty$. Then \eqref{**E} holds for $\kk:= \ff p\dd\ge p,$ so that assertions in Theorem \ref{T1.1} and Theorem \ref{T1.2}(1) hold for $r=p-1.$ When $p=2\ ($i.e. $r:= p-1=1)$ and   $K:=\pi^2-\ff 1 2c_0^2>0,$ Theorem $\ref{T1.2}(2)$ holds.

\subsection{Stochastic generalized fast-diffusion equations}

Let $(E,\B,\m), (L,\D(L)), \H$ and $W_t^{(1)}, W_t^{(2)}$ be in \S 6.1. Let $Q\in \L_{HS}(\H)$ and $B: [0,\infty)\times \V\to \L_{HS}(\H)$ be measurable such that \eqref{LIP} holds for some constant $c_0>0.$

Next, let $r\in (0,1)$, and $\Psi: [0,\infty)\times \R\to\R$ be measurable, continuous in the second variable, such that for some constant $\xi>0$,
\begin{eqnarray}
 && 2\big(\Psi(t,s_1)-\Psi(t,s_2)\big)(s_1-s_2)\ge \ff{\xi |s_1-s_2|^2}{(|s_1|\lor |s_2|)^{1-r}},  \ \ \ s_1,s_2\in \R, t\ge
0,\label{CC1}\\
&& s\Psi(t,s)\ge \xi |s|^{r+1},\ \ \sup_{t\in [0,T],s\ge 0}\ff{|\Psi(t,s)|}{1+|s|^r}<\infty,\ \  s \in \R, t\ge
0,\label{CC2} \end{eqnarray} where  $\ff{|s_1-s_2|^2}{(|s_1|\lor |s_2|)^{1-r}}:=0$ for $s_1=s_2=0.$
By the mean-valued
theorem and   $r\in (0,1)$, one has
$$(s_1-s_2)(s_1^r-s_2^r)\ge
r |s_1-s_2|^2(|s_1|\lor |s_2|)^{r-1},\ \ s_1,s_2\in\R,$$ where $s^r:= |s|^r {\rm sgn}(s)$. So, a  simple example of $\Psi$ for
\eqref{CC1} and \eqref{CC2}  to hold is  $\Psi(t,s)= cs^r$  for some constant $c>0.$

We consider the equation
\beq\label{FD} \d X(t)= \Big\{L\Psi(t,X(t))+ \bb(t) X(t)\Big\}\d t +B(t,X_t)\d W^{(1)}(t)+Q\d W_t^{(2)},\end{equation} where  $\bb\in C([0,\infty))$.   Let $\V=L^{r+1}(\m)\cap \H$ with $\|v\|_\V:= \|v\|_{1+r} + \|v\|$. Then
 $(A1)$-$(A4)$ hold for (see \cite[Theorem 3.9]{RRW} for a more general result) $$A(t,v):= L \Psi(t,v)+ \bb(t) v,\ \ \ v\in \V.$$

\beg{prp}\label{P6.3} Assume $\eqref{LIP}$, $(\ref{CC1})$ and $(\ref{CC2})$. If  there exist
constants $\kk>0$ and $\eta>0$ such that
\beq\label{2.4.3}
\|u\|_{r+1}^{2}  \|u\|^{\kk-2}\ge \eta \|u\|_{Q}^{\kk},\ \
\ u\in L^{r+1}(\m),\ t\ge 0,
\end{equation}then the assertions in Theorem \ref{T1.3} hold. If moreover $B=0$ and $\kk\ge \ff 4 {1+r}$, then $\eqref{NNB}$ holds.
\end{prp}

\beg{proof} It suffices to prove $(A1'')$. By \eqref{CC1} and \eqref{LIP}, there exists constants $K_1,\theta_1>0$ such that
 \beq\label{M3}\beg{split} &  _{\V^*}\<A(t, u)-A(t, v), u-v\>_\V + \ff 1 2  \|B(t, u)-B(t,v)\|_{HS}^2\\
 &\le K_1\|u-v\|^2-\theta_1 \m\big(|u-v|^2 (|u|\lor |v|)^{r-1}\big),\ \ u,v\in \V.\end{split}\end{equation}
 On the other hand, by H\"older's inequality, we have
\begin{eqnarray} &&\|u-v\|_{r+1}^{r+1}= \m(|u-v|^{r+1}) \le \m\big(|u-v|^2(|u|\lor|v|)^{r-1}\big)^{\ff{r+1}2} \m\big((|u|\lor|v|)^{r+1}\big)^{\ff{1-r}2}\no\\
&&\le 2^{\ff{1-r} 2} \m\big(|u-v|^2(|u|\lor|v|)^{r-1}\big)^{\ff{r+1}2} (\|u\|_{1+r}\lor \|v\|_{1+r})^{\ff{1-r^2}2}.\no\end{eqnarray}
Combining this with (\ref{M3}) and \eqref{2.4.3}, we prove $(A1'')$.   \end{proof}

Below we consider the stochastic fast-diffusion equation where   $\Psi(t,s)=s^r := |s|^r {\rm sgn}(s).$

\begin{cor}\label{CFS} Let $\eqref{LIP}$ hold.  Consider $(\ref{FD})$ for   $\Psi(t,s)=s^r:=|s|^r{\rm sgn}(s)$. Let $(-L,\D(L))$ be a nonnegative definite self-adjoint   operator on $L^2(\m)$ with discrete spectrum $(0<)\ll_1\le\ll_2\cdots\le\ll_n\uparrow\infty$ counting multiplicities.   Let $\{e_n\}_{n\ge 1}$ be the corresponding eigenvectors which consist of an orthonormal basis of $L^2(\m)$.  Assume that $-(-L)^{\ff 1 n}$ is a Dirichlet operator for some $n\in\N$ and  the Nash inequality
\beq\label{NASH} \|f\|_{L^2(\m)}^{2+\ff 4 m}\le -C_2 \m(fLf),\ \ f\in\D(L), \m(|f|)=1\end{equation} holds for some constants $C>0$ and $m\in (0, \ff{2(1+r)}{1-r}).$ Let
 $$Q e_i= q_ie_i,\ \ i\ge 1$$ for some constants $\{q_i\}_{i\ge 1}$ satisfying
\beq\label{EI} \|Q\|_{HS}^2= \sum_{i=1}^\infty q_i^2<\infty.\end{equation}
If there exist constants $\kk\ge 2$ and $\vv\in (0,\ff{(1-r)m}{2(1+r)}) $   such that
\beq\label{SB}\sup_{i\ge 1}  |q_i|^{-1}  \ll_i^{ \ff{\vv-1}\kk}<\infty,\end{equation}
then the assertions in Theorem \ref{T1.3} hold. If moreover $B=0$ and $\kk\ge \ff 4 {1+r}$, then $\eqref{NNB}$ holds.
 \end{cor}

\beg{proof}   Obviously,
(\ref{CC1}) and (\ref{CC2})   hold for
 $\Psi(t,s)=s^r.$    To apply Theorem \ref{T1.1}, it remains to verify (\ref{2.4.3}).  By (\ref{SB}) and Schwartz's inequality we have
 \beg{equation*}\beg{split} \|x\|_Q^\kk &= \Big(\sum_{i\ge 1}  \ff{\m(xe_i)^2}{q_i^2\ll_i}\Big)^{\ff\kk 2}\le \Big(\sum_{i\ge 1} \ff{\m(x e_i)^2}{|q_i|^\kk\ll_i}\Big)\Big(\sum_{i\ge 1} \ff{\m(x e_i)^2}{\ll_i}\Big)^{\ff{\kk -2}{2}}\no\\
 &\le c_1   \Big(\sum_{i\ge 1} \ff{\m(x e_i)^2}{\ll_i^\vv}\Big) |x|^{\kk-2}\end{split}\end{equation*} for some constant $c_1>0.$
 According to the proof of Corollary 3.2 in \cite{LW08}, (\ref{NASH}) for some $m\in  (0, \ff{2\vv(r+1)}{1-r})$ implies
 $$\|x\|_{r+1}^2\ge c\sum_{i\ge 1}\ff{ \m(xe_i)^2}{\ll_i^\vv}$$  for some constant $c>0$. Therefore, (\ref{2.4.3}) holds for constant $\eta>0.$
 Then the proof is finished by Proposition \ref{P6.3}\end{proof}

\paragraph{Example 6.3.}\ Let \eqref{LIP} hold, $\Psi(t,s)=s^r$ for some $r\in (0,1)$,    $L=-(-\DD)^\gg$ for some constant $\gg>0$, and  the Dirichlet Laplacian $\DD$ on a bounded domain in $\R^d$, and let $\m$ be the normalized Lebesgue measure on the domain. Let $Qe_i= c_ii^{-\dd} e_i,\ i\ge 1,$ for some constants $\dd\in (\ff 1 2, \ff{1+3r}{2(1+r)})$ and $\{c_i\}_{i\ge 1}\subset \R$ such that $0<\inf_{i\ge 1} |c_i|\le \sup_{i\ge 1}|c_i|<\infty.$
Then for any $\kk\in (\ff{2\gg(1+r)-d(1-r)}{d\dd(1+r)}, \ff{2\gg}{d\dd})\cap [2,\infty)$, assertions in Theorem \ref{T1.3} hold. If moreover $B=0$ and $\kk\ge \ff 4 {1+r}$, then $\eqref{NNB}$ holds.

\beg{proof} We have $\ll_i\ge c i^{\ff{2\gg} d}, i\ge 1$ for some constant $c>0.$ Since $\dd>\ff 1 2$ and $q_i = {\rm O}(i^{-\dd})$, \eqref{EI} holds. Moreover, since $\kk<\ff{2\gg}{d\dd},$ we have $\vv:= 1-\ff{\kk d\dd}{2\gg}\in (0,1)$ and \eqref{SB} follows from $\ll_i\ge c i^{\ff{2\gg}d}. $ Finally, by the classical Nash inequality on $\R^d$ and \cite[Theorem 1.3]{BM}, \eqref{NASH} holds for $m:= \ff d\gg$. Since $\kk> \ff{2\gg(1+r)-d(1-r)}{d\dd (1+r)},$
we have $\vv:= \ff{2\gg- \kk d\dd}{2\gg}\in (0, \ff{(1-r)m}{2(1+r)})$ as required by Corollary \ref{CFS}. Then the proof is finished.\end{proof}

\beg{thebibliography}{99}

\bibitem{ATW} M. Arnaudon, A. Thalmaier, F.-Y. Wang, \emph{Equivalent log-Harnack and gradient for point-wise curvature lower bound,}   Bull. Math. Sci. 138(2014), 643--655.

\bibitem{BBDR} V. Barbu, V.I. Bogachev,  G. Da Prato,  M. R\"ockner, \emph{Weak
solution to the stochastic porous medium  equations: the
degenerate case,}  J. Funct. Anal. 237 (2006), 54--75.

 \bibitem{BDR}  V. Barbu, G. Da Prato, M. R\"ockner, \emph{Finite time extinction of solutions to fast diffusion equations driven by linear multiplicative noise,} J. Math. Anal. Appl. 389(2012), 147--164.

 \bibitem{BM}A. Bendikov, P. Maheux, \emph{Nash type inequalities for fractional powers of non-negative self-adjoint operators,} Trans. Amer. Math. Soc. 359(2007), 3085--3097.

 \bibitem{RW13} M. R\"ockner, F.-Y. Wang, \emph{ General extinction results for stochastic partial differential equations and applications, }   J. Lond. Math. Soc. 87 (2013),   545--560.

\bibitem{CL} M.-F. Chen, S.-F. Li, \emph{Coupling methods for multi-dimensional diffusion processes,} Ann. Probab. 17(1989), 151--177.

\bibitem{DR} G. Da Prato and M. R\"ockner, \emph{Weak solutions to stochastic
porous media equations,} J. Evolution Equ. 4(2004), 249--271.

\bibitem{DRRW} G. Da Prato, M. R\"ockner, B.L. Rozovskii, F.-Y. Wang, \emph{Strong solutions of Generalized porous media equations: existence, uniqueness and ergodicity, } Comm. Part. Diff. Equ. 31 (2006), 277--291.
\bibitem{Gess} B. Gess, \emph{Finite time extinction for stochastic sign fast diffusion and self-organized criticality,} arXiv: 1310.6971.

\bibitem{KR}
  N.V. Krylov, B.L. Rozovskii,
  \emph{Stochastic evolution equations,}
  Translated from Itogi Naukii Tekhniki, Seriya Sovremennye Problemy Matematiki,
  14(1979), 71--146, Plenum Publishing Corp. 1981.

 \bibitem{LR} T. Lindervall, L.C.G. Rogers, \emph{Coupling of multidimensional diffusions by reflection,} Ann. Probab. 14(1986), 860--872.

\bibitem{L} W. Liu, \emph{Harnack inequality and applications for stochastic evolution equations with monotone drifts,}  J. Evol. Equ. 9(2009),  747--770.

\bibitem{LR1} W. Liu, M. R\"ockner, \emph{SPDE in Hilbert space with locally monotone coefficients,} J. Differential Equations 259(2010), 2902--2922.

\bibitem{LR2} W. Liu, M. R\"ockner, \emph{Local and global well-posedness of SPDE with generalized coercivity conditions,} J. Funct. Anal. 254(2013), 725--755.

\bibitem{LW08} W. Liu, F.-Y. Wang, \emph{Harnack inequality and strong Feller property for stochastic fast-diffusion equations,} J. Math. Anal. Appl. 342(2008), 651--662.

\bibitem{P1} E. Pardoux, \emph{Sur des equations aux d\'eriv\'ees
partielles stochastiques monotones,} C. R. Acad. Sci. 275(1972),
A101--A103.

\bibitem{P2} E. Pardoux, \emph{Equations aux d\'eriv\'ees
partielles stochastiques non lineaires monotones: Etude de solutions
fortes de type Ito,} Th\'ese Doct. Sci. Math. Univ. Paris Sud. 1975.

\bibitem{PW} E. Priola, F.-Y. Wang, \emph{Gradient estimates for diffusion semigroups with singular coefficients,}  J. Funct. Anal. 236(2006), 244--264.

\bibitem{RRW} J. Ren, M. R\"ockner, F.-Y. Wang, \emph{Stochastic generalized porous media and fast diffusion equations,} J. Differential Equations 238(2007), 118--152.

\bibitem{RW10} M. R\"ockner, F.-Y. Wang, \emph{Harnack and functional inequalities for generalized Mehler semigroups,}  J. Funct. Anal.  203(2007), 237--261.

\bibitem{W07} 	F.-Y. Wang, \emph{Harnack inequality and applications for stochastic generalized porous media equations,}  Annals of Probability 35(2007), 1333--1350.

\bibitem{Wbook} F.-Y. Wang, \emph{Harnack Inequality and Applications for Stochastic Partial Differential Equations,} Springer, New York, 2013.

\bibitem{W14} F.-Y.  Wang, \emph{ Exponential convergence of non-linear monotone SPDES, } arXiv: 1310.7997v2.

\bibitem{WZ} F.-Y. Wang, T. Zhang,  \emph{Log-Harnack inequalities for semi-linear SPDE with strongly multiplicative noise,}  Stoch. Proc. Appl. 124(2014), 1261--1274.

\end{thebibliography}

\end{document}